\documentclass[11pt,a4paper]{article}
\usepackage{a4wide}
\usepackage{mathrsfs} 
\usepackage{amssymb} 
\usepackage{amsmath} 
\usepackage{amsthm} 
\usepackage{hyperref}

\usepackage{enumitem}
\setlist[itemize]{leftmargin=*}
\setlist[enumerate]{leftmargin=*}

\makeatletter
\renewcommand\@biblabel[1]{#1.}
\makeatother

\usepackage{fancyhdr}
\date{\vspace{-\baselineskip}}

\makeatletter
\renewcommand\tableofcontents{%
    \@starttoc{toc}%
}
\makeatother
\usepackage{titletoc}%
\contentsmargin{0pt}
\titlecontents*{section}[0pt]
{\begin{quote}\small}{\thecontentslabel. }{\end{quote}}
{\\}[][References]

\usepackage[center]{titlesec}%

\titleformat{name=\section}
{\Large\center\bfseries}{\thesection.}{1em}{}
\titleformat{name=\subsection}
{\large\center\bfseries}{}{}{}
\titleformat{name=\subsubsection}
{\large\center\em}{}{}{}

\newtheoremstyle{dotless}{}{}{\itshape}{}{\bfseries}{}{ }{}
\theoremstyle{dotless}

\newtheorem{theorem}{Theorem}[section] 
\newtheorem{lemma}[theorem]{Theorem} 
\newtheorem*{corollary}{Corollary}
\newtheorem{conjecture}{Conjecture}

\newtheorem*{LLL}{Lov\'asz local lemma}
\newtheorem*{definition}{Key Definition}
\newtheorem*{question}{Question}

\usepackage{tcolorbox}
\usepackage{algorithm}
\usepackage{algorithmic}
\usepackage{mathrsfs}

\newcommand*{\EE}{\mathbb{E}}
\newcommand*{\cI}{\mathcal{I}}

\renewcommand{\epsilon}{\varepsilon}
\newcommand{\eps}{\varepsilon}
\renewcommand{\Pr}{\mathbb{P}}
\renewcommand{\subset}{\subseteq}

\newcommand{\avdeg}{\rho}

\newcommand{\tomchapter}{7}
\newcommand{\ourchapter}{8}

\title{{\bf \ourchapter}\\The hard-core model in graph theory
\footnote{This survey has been prepared as a chapter for a forthcoming volume, \emph{Topics in Probabilistic Graph Theory}. Draft dated 3 September 2025.}
}

\author{
{\sc Ewan Davies
\footnote{Department of Computer Science, Colorado State University, Fort Collins, USA. \protect\href{mailto:research@ewandavies.org}{\protect\nolinkurl{research@ewandavies.org}}}
}
\ \ and \ \
{\sc Ross J. Kang
\footnote{Korteweg--de Vries Institute for Mathematics, University of Amsterdam, Netherlands. \protect\href{mailto:ross.kang@gmail.com}{\protect\nolinkurl{ross.kang@gmail.com}}}
} 
}

\begin{document}

\maketitle

\tableofcontents

\begin{abstract}
\noindent\em
An independent set may not contain both a vertex and one of its neighbours. This basic fact makes the uniform distribution over independent sets rather special. We consider the hard-core model, an essential generalization of the uniform distribution over independent sets. We show how its local analysis yields remarkable insights into the global structure of independent sets in the host graph, in connection with, for instance, Ramsey numbers, graph colourings, and sphere packings.
\end{abstract}



\section{Introduction}\label{sec:intro}

For the uninitiated, the hard-core model might sound \ldots hard-core, so to speak.

In this chapter, our aim is not to give the most comprehensive or in-depth treatment of the use of the hard-core model in graph theory, which would indeed be a bit hard-core (and perhaps self-defeating), but rather to give a current overview of some of its most elegant and interesting uses, and to provide a glimpse into its varied applicability. In this section, with the help of a `toy' example, we introduce the model. Then in the following section, we sketch a series of `vignettes', where the hard-core model, and in particular the local occupancy method, has had influence. In later sections, we dive into more of the details of a recent framework for inferring graph structure that was developed around the hard-core model. At the end, we discuss some frontiers related to the methods and topics covered in this chapter.


\subsection*{A toy example}

Given a graph, consider a uniformly random independent set $X$. To be completely clear, we first gather all of the independent sets of the graph (including $\varnothing$) into a collection $\cI$, and then choose $X$ randomly from  $\cI$, where each member of  $\cI$ is chosen with equal probability.
Let $u$ be a vertex and reveal the number  $Y_u=|X\cap N(u)|$ of neighbours of $u$ {\em occupied} by $X$. 
Since $X$ is an independent set, the event $u\in X$, that $u$ itself is occupied, cannot occur unless $Y_u=0$.
Moreover, since $X$ is chosen uniformly from $\cI$, it is intuitive that, conditional on $Y_u=0$, the events $u\in X$ and $u\notin X$ have equal probability (and we prove this rigorously in more generality below).
Thus, $\Pr(u\in X)= \frac{1}{2}\Pr(Y_u=0)=\tfrac{1}{2}q$, where we write $q=\Pr(Y_u=0)$.
Markov's inequality then implies that $\EE|X\cap N(u)|=\EE Y_u\ge \Pr(Y_u\ge 1) = 1-q$. It follows that 
\[ 2\Pr(u\in X) + \EE|X\cap N(u)|\ge 2\cdot\tfrac{1}{2}q +(1-q) = 1. \]
The two quantities on the left-hand side are in tension. Since $X$ is an independent set, one quantity is the probability that $u$ is occupied, and the other is the sum of probabilities of the neighbours of $u$ being occupied.
As such, bounding from below a positive linear combination of these two quantities can be seen as `balancing' them.

This `balancing' (uniformly over all $u$) immediately gives a lower bound on the average size of an independent set.
Writing $n$ for the number of vertices and $\Delta$ for the maximum vertex-degree, we sum the inequality over all vertices $u$ to obtain
\begin{align*}
    n\le 2\sum_{u}\Pr(u\in X) + \sum_u\sum_{v\sim u}\Pr(v\in X) \le (\Delta+2)\EE|X|,
\end{align*}
because $\EE|X|=\sum_u\Pr(u\in X)$ and in the double sum each vertex $u$ appears at most $\Delta$ times. 
(Here $v\sim u$ denotes that the vertices $v$ and $u$ are adjacent.)
Rewriting this, we have
\begin{align*}
\EE|X| \ge \frac{n}{\Delta+2}.
\end{align*}
It might not be obvious to see how good this bound is, but we will see that it is tight.
Note that, by the probabilistic method, there is guaranteed to be an independent set of size at least the quantity on the right-hand side.

Let us now illustrate a generalization of this with an appropriate `rescaling' of the distribution of $X$ according to some parameter $\lambda>0$. 
Consider a random independent set $X$ of the graph, where, instead of a uniform choice, any $I\in \cI$ is chosen with probability proportional to $\lambda^{|I|}$ -- that is, the probability of choosing $I \in \cI$ is
\begin{align*}
\Pr(X = I) = \frac{\lambda^{|I|}}{\sum_{J\in\cI} \lambda^{|J|}}.
\end{align*}
This distribution over the collection $\cI$ of all independent sets is also known as the {\em hard-core model at fugacity $\lambda$}.
The uniform case corresponds to the choice of fugacity $\lambda=1$. 
The normalizing factor
\[
Z(\lambda) = \sum_{J\in\cI} \lambda^{|J|}
\]
is known as the {\em partition function} of the hard-core model or the {\em independence polynomial} of the graph, and it is an important parameter of the system.

As before, let $u$ be any vertex and reveal $Y_u=|X\cap N(u)|$. As $X$ is an independent set, the event that $u\in X$ may not occur unless $Y_u=0$. It turns out, as we shall shortly discuss in more detail, the `rescaling' of the uniform distribution on $\cI$ that we chose enjoys something remarkable called the `spatial Markov property'. This property implies that, conditional on the event $Y_u=0$, the probability of the vertex $u$ being chosen in $X$ is the same as that of the one-vertex independent set being chosen in a one-vertex graph under the same rescaling -- that is, with probability $\lambda/(1+\lambda)$. (This is just as we observed above in the case $\lambda=1$.)
Writing $q=\Pr(Y_u=0)$, we then have $\Pr(u\in X)= \lambda q/(1+\lambda)$.
Again, we have  $\EE|X\cap N(u)|\ge 1-q$ by Markov's inequality, and so 
\begin{align} 
\frac{1+\lambda}{\lambda}\,\Pr(u\in X) + \EE|X\cap N(u)|\ge 1, \label{eq:lococcgreedy}
\end{align}
a `balancing' which we see later as a consequence of a property that we call \emph{local $((1+\lambda)/\lambda,1)$-occupancy}. 
Summing this inequality over all vertices $u$ and using the fact that \[\sum_u\sum_{v\sim u}\Pr(v\in X) \le \Delta \cdot \EE|X|\] yields
\begin{align}
n \le (\Delta+1+1/\lambda)\cdot \EE|X|
\quad\text{ or, equivalently, }\quad
\EE|X| \ge \frac{n}{\Delta+1+1/\lambda}. \label{eq:occgreedy}
\end{align}
Again, by the probabilistic method, there is guaranteed to be an independent set of size at least the quantity on the right-hand side. 
We generalized what we did earlier, but now we are free to choose $\lambda$ as large as we like. Indeed it is valid to take the limit $\lambda\to\infty$, and so there must be an independent set of size at least $n/(\Delta+1)$.

It seems that up to here we have expended substantial energy to prove a simple independence number bound in bounded-degree graphs. Note that, by the greedy method of arbitrarily adding a vertex to $I$ and then deleting it and all its neighbours from the graph, we can easily find an independent set $I$ of size at least  $n/(\Delta+1)$.
Indeed, {\em every} maximal independent set has size at least this value.

 However, we point out how~\eqref{eq:occgreedy}, as stated, is actually best possible. To see why, consider a complete graph on $n=\Delta+1$ vertices. The collection $\cI$ then comprises the empty set (of relative probability mass $1$) and $\Delta+1$ independent sets of size $1$ (each of relative probability mass $\lambda$), and so the expected size of $X$ in this example satisfies
\begin{align*}
\EE|X| = \frac{1\cdot0+(\Delta+1)\lambda}{1\cdot\lambda^0+(\Delta+1)\lambda^1}=\frac{\Delta+1}{\Delta+1+1/\lambda}=\frac{n}{\Delta+1+1/\lambda}.
\end{align*}

We have taken this longer route in the basic setting of general bounded-degree graphs in an effort to hint at the value an inequality of the type in~\eqref{eq:lococcgreedy} can have in deriving good bounds on the independence number.
Later applications will show how we can carefully leverage local graph structure and optimize other bounds of the form in~\eqref{eq:lococcgreedy} with respect to the hard-core model, and in many cases produce state-of-the-art bounds. 
We shall demonstrate, moreover, that this type of local occupancy bound, especially in concert with the spatial Markov property, is sufficient for finding nice and `more even' distributions over $\cI$, such as good fractional and integral proper colourings. 

\subsection*{Additional context}

To provide more motivation and additional context for the hard-core model, let us see how the above `scaling' could be derived naturally from some simple probabilistic considerations.
Let $A\subset V(G)$ be a random vertex-subset of $G$, formed by including each vertex of $G$ independently at random with probability $p \in(0,1)$.
Now let $X$ be $A$ conditioned on the event that $A$ forms an independent set of $G$. Given $I\in \cI$, we compute
\[ \Pr(X=I) = \Pr(A = I \mid A \in\mathcal{I}) = \frac{\Pr(A=I)}{\Pr(A\in \mathcal{I})} = \frac{(1-p)^{|V(G)|}}{\Pr(A\in\mathcal I)} \cdot \left(\frac{p}{1-p}\right)^{|I|},\]
and so the only information about $I$ on which this probability depends is $|I|$. We have recovered the hard-core model on $G$ at fugacity $p/(1-p)$.

Further motivation comes from the fact that it is a fundamental example of a \emph{Gibbs measure}. 
Specifically, the distribution on $\mathcal I$ which maximizes entropy, subject to a fixed expected size, turns out to be the hard-core model. 
Given $I\in \cI$, let $p_I = \Pr(X=I)$. The information entropy of the random variable $X$ is defined as $H(X) = -\sum_{I\in\mathcal I} p_I \log_2 p_I$. 
The method of Lagrange multipliers tells us that the distribution which fixes $\mathbb{E}|X| =\mu$ and maximizes entropy is a stationary point of the function
\[ \mathcal{L} = -\sum_{I\in\mathcal I}p_I\log_2 p_I  + \Lambda_1\left(\sum_{I\in\mathcal I}p_I -1\right) +\Lambda_2\left(\sum_{I\in\mathcal I}|I|p_I - \mu\right).\]
A little calculus shows then that $p_I = \exp(-1 + \Lambda_1\log 2 + |I|\Lambda_2 \log 2)$, and again there must be some value of $\lambda$ such that $p_I\propto \lambda^{|I|}$.
(Unless specified otherwise, the base of the logarithm is always natural.)
The precise value of $\lambda$ is a function of $\mu$, and for the rest of the calculation, we require that $0\le \mu\le \alpha(G)$. 
These properties mean that the hard-core model arises naturally in physics as a discrete model of gas particles that are constrained to the vertices of $G$. 
The only interaction that these particles have is that they cannot simultaneously occupy both endpoints of an edge of $G$, and this self-exclusion property forms the basis of a discrete model of a gas (the so-called \emph{hard-core lattice gas}) where particles have fixed size and cannot overlap, but are not subject to inter-particle forces. For more of a physics perspective, see, for example,~\cite{BS94},~\cite{FV17}, and~\cite{Geo11}.

\section{Vignettes}\label{vignettes}

In this section, we sketch a series of `vignettes', where the hard-core model, and in particular the local occupancy method, has had influence.

\subsection*{Off-diagonal Ramsey numbers}

A natural starting vignette is the following classical combinatorial problem. The {\em off-diagonal Ramsey number} $R(3,k)$ is the least integer $n$ for which, in any $2$-edge-coloured complete graph on $n$ vertices, there is guaranteed to be either a triangle with all of its edges in the first colour or a complete subgraph on $k$ vertices with all of its edges in the second colour.
Stated in another way, $R(3,k)$ is the least integer $n$ for which any triangle-free graph on $n$ vertices must have an independent set of size $k$.
This notion goes back a century to the work of Ramsey~\cite{ramsey_problem_1930}, who required the well-definedness of some more general parameters (the multicolour hypergraph Ramsey numbers) in order to show the decidability of some statements in first-order logic. The quest for better quantitative bounds for this and related notions has deeply influenced the development of combinatorics in general, and especially probabilistic and extremal combinatorics (see~\cite{Erd61},~\cite{Spe77}, and~\cite{Kim95}). 
There have been relatively recent breakthroughs, through the sophisticated analysis of the so-called `triangle-free process', for lower bounds on $R(3,k)$ (see~\cite{Boh},~\cite{BoKe}, and~\cite{FGM}). 
On the other hand, an important demonstration of hard-core methods matches the best-known {\em upper} bounds for $R(3,k)$.

It is simple, but crucial, to notice that in any triangle-free graph, the set of neighbours of any given vertex forms an independent set. 
It is therefore intuitive that triangle-freeness might lead to there being larger independent sets than in general, as there are independent sets spread locally throughout the graph!
Indeed, if we are interested in finding an independent set of size $k$ for any triangle-free graph on $n$ vertices, we may at least assume without loss of generality that all vertices in the graph have degree less than $k$.
We have thereby reduced the problem of showing $R(3,k)\le n$ to that of finding an independent set of size $k$ in any triangle-free graph on $n$ vertices of maximum degree at most $k-1$.
It turns out that the more general problem of finding large independent sets in triangle-free graphs on $n$ vertices of maximum degree at most $\Delta$ is the crux of bounding $R(3,k)$ from above, especially asymptotically, as $k\to\infty$.

\begin{question}
In a triangle-free graph on $n$ vertices of maximum degree at most $\Delta$, how large an independent set can we guarantee?
\end{question}


The seminal work of Ajtai, Koml\'os and Szemer\'edi~\cite{AKS80},~\cite{AKS81} -- which simultaneously introduced an early form of the nibble method -- showed in this setting that the independence number must be at least $\Omega((n\log\Delta)/\Delta)$, improving on the `trivial' bound presented in Section~\ref{sec:intro} by a factor logarithmic in $\Delta$. 
Not long afterwards, Shearer~\cite{She83} improved this independence number bound by a constant factor to $(1+o(1))(n\log\Delta)/\Delta$ 
as $\Delta\to\infty$.
If we substitute $\Delta=k-1$ as in the previous paragraph, then we should clearly take $n=(1+o(1))k^2/\log k$ to derive that $R(3,k) \lesssim k^2/\log k$. 
Note that these independence number results apply in the more general setting of given \emph{average} degree, though one does not need this additional strength for the Ramsey number bound.

Despite sustained attention from many leading researchers, Shearer's bound remains asymptotically the state-of-the-art upper bound for $R(3,k)$.
The earlier-mentioned work on the triangle-free process (see~\cite{BoKe},~\cite{FGM}) shows this to be counterbalanced by a lower bound of the form $R(3,k) \gtrsim k^2/(4\log k)$.

We cannot fail to mention how special the case of forbidding {\em triangles} seems to be. As one illustration, Ajtai, Erd\H{o}s, Koml\'os and Szemer\'edi~\cite{AEKS81} conjectured that every $K_4$-free graph of $n$ vertices of maximum degree at most $\Delta$ must contain an independent set of size $\Omega((n\log\Delta)/\Delta)$; however, the best bound is of size $\Omega((n\log\Delta)/(\Delta\log\log\Delta))$ (see~\cite{She95}), with the current best asymptotic leading constant in~\cite{davies_graph_2020}.

One of the examples of local occupancy analysis we present gives an alternative derivation, first published in~\cite{DJPR17average},
 of Shearer's $R(3,k)$ bound.

\subsection*{Sphere packing}

The sphere packing problem is to determine, for each dimension $d$, the density $\theta(d)$ of the densest possible packing of congruent spheres in $\mathbb{R}^d$. 
This is clearly a natural geometric problem, and constitutes a basic question about models of granular materials. 
The case of dimension $d=1$ is trivial, and of $d=2$ elementary, but it seems that every other case is challenging. 
The optimality of the familiar packing of 3-dimensional spheres that one might use to stack oranges was conjectured by Kepler~\cite{kepler_strena_1611} in 1611 and was proved in 2005 by Hales~\cite{hales_proof_2005},~\cite{hales_formal_2017}. 
The dimensions $d=8$ and $d=24$ required another breakthough due to Viazovska~\cite{viazovska_sphere_2017},~\cite{cohn_sphere_2017}, but $\theta(d)$ is not known exactly for any other dimension.
The handful of packings that are known to be densest arise from lattices, which means that the centres of the spheres are arranged regularly throughout the space $\mathbb{R}^d$.

\begin{question}
Given a collection of non-intersecting congruent spheres in $\mathbb{R}^d$, how high can their density in $\mathbb{R}^d$ be?
\end{question}


The importance of high-dimensional sphere packings, beyond natural geometric curiosity, became apparent through the work of Shannon~\cite{shannon_mathematical_1948} on coding theory. 
The centres of a packing of spheres correspond to an error-correcting code for a simple model of communication through a noisy channel, and high-density packings yield codes with better information-carrying capacity.
It is easy to prove the so-called `Minkowski bound', that a packing of density at least $2^{-d}$ is possible, since any maximal packing must have the property that if we double the radii of the spheres, then all of $\mathbb{R}^d$ is covered. 
From a graph-theoretic perspective, we can choose a radius $r$ and consider a graph on $\mathbb{R}^d$, where the vertices $u$ and $v$ are adjacent if and only if $d(u,v)\le 2r$. 
Independent sets in this graph correspond to packings of non-overlapping spheres of radius $r$, and analyzing the density of such packings is formidable, but tractable. 

Krivelevich, Litsyn and Vardy~\cite{krivelevich_lower_2004} studied a discrete subgraph of this infinite graph to match the growth rate $\Omega(d2^{-d})$ of the best lower bounds on $\theta(d)$ known at the time. 
They discretized the problem by taking the vertices on the integer lattice intersected with a finite cube, and bounded the independence number of the resulting graph in terms of its maximum degree. 
Using the trivial bound $\alpha(G)\ge |V(G)|/(\Delta+1)$ recovers the Minkowski bound on $\theta(d)$, but unless the graph is especially dense, this bound is far from tight. 
One might hope that the graph is triangle-free, which by Shearer's result mentioned above would yield an improvement by a factor $\Omega(\log\Delta)$ and by a corresponding factor $\Omega(d)$ in the lower bound on $\theta(d)$, but it seems that taking the vertices suitably close together in $\mathbb{R}^d$ necessitates some triangles. 
Fortunately, the graphs considered contain few enough triangles that the desired improvements come from generalizations of the triangle-free independence number bounds to so-called `locally sparse' graphs, which can be proved with the hard-core model. 

Inspired by the corresponding finite results, Jenssen, Joos and Perkins~\cite{jenssen_hard_2019} studied directly the analogue of the hard-core model on the infinite graph on $\mathbb{R}^d$, known as the \emph{hard-sphere model}. They proved that a random packing from this model achieves a density of $\Omega(d2^{-d})$ in expectation, and
their techniques evoked the concept of local occupancy that we study here in the continuous setting.
The best-known lower bound $\theta(d)\ge \Omega(d2^{-d}\log d)$ is due to Campos, Jenssen, Michelen and Sahasrabudhe~\cite{campos_new_2023} and relies on a novel disordered discretization of space, as well as a new lower bound on the independence number of locally sparse graphs.
Apart from determining $\theta(d)$, the question of whether dense packings in high dimensions are disordered or lattice-like remains elusive.

\subsection*{Chromatic number of triangle-free graphs}

By a simple greedy algorithm similar to the one for independent sets, where we consider vertices one at a time and assign the least available colour, we see that the chromatic number of a graph of maximum degree $\Delta$ is at most $\Delta+1$.
A classic result of Brooks~\cite{brooks_colouring_1941} states that this is attained only for graphs containing cliques on $\Delta+1$ vertices or odd cycles when $\Delta=2$. Another way of seeing Brooks's theorem is that, for $\Delta\ge 3$, we strictly improve upon the greedy bound upon the exclusion of a clique of size more than $\Delta$. In 1968, Vizing~\cite{vizing_unsolved_1968} asked about best chromatic number bounds upon the exclusion of smaller cliques, and triangles, especially.

\begin{question}
Given a triangle-free graph of maximum degree $\Delta$, how large can the chromatic number be?
\end{question}


In general, this problem is very difficult. Even the case $\Delta=5$ remains open:  in 1970, Gr\"unbaum~\cite{Gru70} conjectured that there exists some triangle-free $5$-regular graph of chromatic number $5$. While others conjectured the opposite (see~\cite{Ree98}), it remains open to debate. For other small values of $\Delta$, an upper bound of $2\lceil (\Delta+2)/3\rceil$ is due to Kostochka in 1982 (see~\cite[Sec.~4.6]{jensen_graph_1994}).

In a remarkable unpublished work in the 1990s, using the nibble method (see Chapter~\tomchapter), Johansson~\cite{johansson_asymptotic_1996} (see also~\cite{molloy_graph_2002}) showed how in triangle-free graphs the asymptotic extremal behaviour,  for large $\Delta$, of the chromatic number is aligned to that of the independence number, and showed that the chromatic number may be no larger than $O(\Delta/\log\Delta)$. This implies the Ajtai--Koml\'os--Szemer\'edi result by taking a largest colour class in an optimal proper colouring.

In a recent conceptual advance, Molloy~\cite{molloy_list_2019} incorporated the use of the entropy compression method due to Moser~\cite{Mos09} into a local randomized recolouring procedure for triangle-free graphs. Not only did he achieve a better bound, of the form $(1+o(1))\Delta/\log \Delta$, but he also discovered a substantially shorter and more elegant proof than Johansson's. Again, by taking a largest colour class, this chromatic number bound matches the longstanding bound set by Shearer for the independence number, which is at the same time a longstanding bound for the off-diagonal Ramsey numbers $R(3,k)$.

In Molloy's proof, the local use of a uniformly random partial list colouring is critical, and in fact relies on the spatial Markov property! This is also critical to the local occupancy framework that we outline in this chapter.

\subsection*{Counting independent sets}

At a number theory conference in Banff in 1988, Granville asked for an upper bound on the number of independent sets in an $n$-vertex $d$-regular graph -- specifically one of the form $2^{(1/2 +o_d(1))n}$, where the $o_d(1)$ term is as $d\to 0$.
Granville was motivated by a problem due to Cameron~\cite{cameron_portrait_1987}, who conjectured that the number of subsets of $\{1,2,\dotsc,n\}$ which contain no solutions to $x+y=z$ is at most $2^{(1/2+o(1))n}$. 
We say that the set $A$ is \emph{sum-free} if there are no solutions in $A$ to the equation $x+y=z$. 

The connection between the problems is somewhat subtle. 
For an Abelian group $\Gamma$ of size $n$ and a subset $S\subset \Gamma$, we say that $A\subset \Gamma$ is \emph{$S$-free} if the set $A+S = \{a+s : a\in A, s\in S\}$ has an empty intersection with $A$.
The $S$-free subsets of $\Gamma$ correspond to independent sets in the \emph{Cayley graph} $G_S$ with vertex-set $\Gamma$ such that $x$ and $y$ are adjacent if and only if either $x-y \in S$ or $y-x\in S$. 
Such Cayley graphs are $d$-regular with $d=|S\cup (-S)|$.
Unless we use addition modulo $n$, the set $[n]$ is not an Abelian group under addition. This technicality can be worked around, but we omit the details.

When $|\Gamma|=n$, there are at most 
$2^{o(n)}$ subsets of $\Gamma$ of size less than $\log_2 n$, and so it suffices to obtain the bound $2^{(1/2 +o(1))n}$ on the number of sum-free subsets of $\Gamma$ of size at least $\log_2n$.
For each such sum-free subset $A$ of $\Gamma$, let $S_A$ be the first $\log_2 n$ elements of $A$. 
Then $A$ is an $S_A$-free subset of $\Gamma$, and so for each fixed $S$ we want an upper bound on the number of independent sets in $G_S$, which is regular of degree $|S\cup (-S)|\ge \log_2 n$. 
Supposing that we have established Granville's conjectured bound, we have that each $S$ can arise as the initial segment of at most $2^{(1/2 +o(1))n}$ sum-free sets $A$. 
The number of possible sets $S$ is $2^{o(n)}$, giving the bound $2^{(1/2 + o(1))n}$ on the number of large sum-free sets, as required.

\begin{question}
Given a $d$-regular graph on $n$ vertices, how large can the number of independent sets be?
\end{question}


Several solutions to Cameron's problem appeared around 1990 (see~\cite{calkin_number_1990},~\cite{alon_independent_1991} and an unpublished proof due to Erd\H{o}s and Granville~\cite{CE90}). But it was Alon~\cite{alon_independent_1991} who tackled Granville's problem on independent sets in regular graphs. 
Alon suggested that a sharper bound is plausible: that when $n$ is divisible by $2d$, the tight upper bound, which is attained by a disjoint union of $n/(2d)$ copies of the complete bipartite graph $K_{d,d}$, should be $(2^{d+1}-1)^{n/(2d)}$. 
With a different motivation, related to statistical physics and phase transitions of the hard-core model, Kahn~\cite{kahn_entropy_2001} developed the so-called `entropy method' and established this bound in the case of bipartite regular graphs. 
Later, Zhao~\cite{zhao_number_2010} reduced the general case to the bipartite case, thereby settling the problem for all regular graphs. 
These results were extended to bounds on the partition function of the hard-core model and to more general weighted counting problems by Kahn~\cite{kahn_entropy_2002} and Galvin and Tetali~\cite{galvin_weighted_2004}. 

For the hard-core model, the results above were strengthened to bounds on the expected size of the independent set by Davies, Jenssen, Perkins and Roberts~\cite{DJPR17independent}.
Their method readily incorporated local structure, which offers advantages in some settings. 
Their methods could be applied to triangle-free graphs (see~\cite{DJPR17average}) to produce independence number bounds analogous to Shearer's. In fact, Shearer knew that the hard-core model could be used in this way and presented his work at the 1998 SIAM annual meeting and conference on discrete mathematics held in Toronto~\cite{She98}.
One of the open problems that Shearer highlighted was the question of whether the analysis of the hard-core model could provide an efficient algorithm for constructing independent sets of the size certified by the method and, through the work of Molloy~\cite{molloy_list_2019} and Bonamy, Kelly, Nelson and Postle~\cite{bonamy_bounding_2022}, these threads were finally united in the form of local occupancy (see~\cite{DJKP20loc}).

\section{Random independent sets: a balancing act}\label{sec:balancing}

Our objective in this section is to motivate how various interesting local properties of the hard-core model at fugacity $\lambda$ lead naturally to nice bounds on the occupancy fraction, the independence number and the fractional chromatic number, subject to an analytic optimization in terms of $\lambda$ and the maximum degree.
We remark that, while the power to optimize over $\lambda$ is handy, the specialization to the uniform case $\lambda=1$ already contains much of the important intuition.
We have given a hint of this line of reasoning in Section~\ref{sec:intro}, but we now provide extra detail and intuition.

In this, we draw some important approaches close together, including those of Shearer~\cite{She95}, Alon~\cite{alon_independence_1996}, Molloy and Reed~\cite{molloy_graph_2002}, and a few later works involving the authors (\cite{DJPR17average}, \cite{DJKP20loc}, \cite{davies_occupancy_2021}, \cite{davies_graph_2020}).

Here follow some standing assumptions for this section. 

\begin{itemize}
\item
$G=(V,E)$ is a graph of maximum degree $\Delta$.
\item
$\cI=\cI(G)$ is the collection of independent sets in $G$ (including $\varnothing$).
\item
For a fixed $\lambda>0$,
$X$ is chosen randomly from $\cI$ according to the following distribution: for any $I\in \cI$,
\vspace{-6pt}
\begin{align*}
\Pr(X = I) = \frac{\lambda^{|I|}}{\sum_{J\in\cI} \lambda^{|J|}}.
\end{align*}
\end{itemize}

%
%

\noindent 
As we mentioned earlier, this distribution over $\cI$ is known as the {\em hard-core model on $G$ at fugacity $\lambda$},
and we refer to the normalizing factor
\[
Z_G = Z_G(\lambda) = \sum_{J\in\cI} \lambda^{|J|}
\]
as its {\em partition function}. 

It is meaningful to observe that the derivative of this function satisfies
\begin{equation}\label{eq:Z'}
\lambda Z_G'(\lambda) = {\sum_{J\in\cI} |J| \lambda^{|J|}} \quad\text{ and so }\quad \EE |X| = \frac{\lambda Z_G'(\lambda) }{Z_G(\lambda)}.
\end{equation}

Many of the methods that we describe here take advantage of an elegant Markovian property of the hard-core model, which allows us to isolate our considerations between different parts of the graph.
Suppose that we want to know how $X$ is distributed with respect to the subgraph $G[A]$ induced by some subset $A\subset V$.
Of course, if we expose $X$ outside of $A$, then we must exclude consideration of those  vertices in $A$ that have some neighbour (external to $A$) occupied by $X$.
The {\em spatial Markov property} asserts that it is sufficient to independently instantiate the distribution on those vertices in $A$ that remain, the so-called \emph{externally uncovered vertices} in $A$.
That assertion is formalized in the following result.

\begin{lemma}[Spatial Markov property of the hard-core model]\label{lem:spatialMarkov}
For any $A\subset V$, $X\cap A$ is distributed according to the hard-core model on the (random) graph $F_A=G[A\setminus N(X\setminus A)]$ at fugacity $\lambda$.
\end{lemma}
To see this, let $I_0$ be an arbitrary independent set of $G[V\setminus A]$, and let us condition on the event that $X\setminus A=I_0$. 
It follows that $X\cap A$ is a subset of $A\setminus N(I_0)$, and indeed is an independent set in $G[A\setminus N(I_0)]=F_A$. 
For any independent set $I_1$ in $F_A$, we have
\begin{align*}
    \Pr(&X\cap A=I_1 \mid X\setminus A=I_0) 
    \\&= \frac{\Pr(X=I_0\cup I_1 \wedge X\setminus A = I_0)}{\Pr(X\setminus A=I_0)} = \frac{\Pr(X=I_0\cup I_1)}{\Pr(X\setminus A=I_0)}
    \\&= \frac{\lambda^{|I_0\cup I_1|}}{Z_G(\lambda)}\cdot\frac{Z_G(\lambda)}{\sum_{J\in\cI(F_A)}\lambda^{|J\cup I_0|}}
    = \frac{\lambda^{|I_1|}}{\sum_{J\in\cI(F_A)}\lambda^{|J|}},
\end{align*}
and this completes the proof of Theorem~\ref{lem:spatialMarkov} because $Z_{F_A}(\lambda) = \sum_{J\in\cI(F_A)}\lambda^{|J|}$.
\medskip

We can now derive a straightforward corollary.

\begin{corollary}
For any $A\subset V$,  the following hold for $F_A=G[A\setminus N(X\setminus A)]$:
\begin{align}\label{eqn:hardcore}
 \Pr(X\cap A = \varnothing) = \EE\frac{1}{Z_{F_A}(\lambda)}\quad\text{ and }\quad
\EE|X\cap A| = \EE\frac{\lambda Z_{F_A}'(\lambda)}{Z_{F_A}(\lambda)}.
\end{align}
\end{corollary}

For this derivation, we let $A_0$ be an arbitrary subset of $A$, and we condition on the event that $A\setminus N(X\setminus A) = A_0$.
By Theorem~\ref{lem:spatialMarkov}, we have
\begin{align*}
    \Pr(X\cap A=\varnothing \mid A\setminus N(X\setminus A) = A_0)
    &= \frac{1}{Z_{G[A_0]}(\lambda)}
\text{ and}\\
    \EE[|X\cap A| \mid A\setminus N(X\setminus A) = A_0]
   & = \frac{\lambda Z_{G[A_0]}'(\lambda)}{Z_{G[A_0]}(\lambda)}.
\end{align*}
It follows that
\begin{align*}
\Pr(X\cap A=\varnothing) &= \sum_{A_0\subset A} \frac{1}{Z_{G[A_0]}(\lambda)} \Pr(A\setminus N(X\setminus A) = A_0) = \EE\frac{1}{Z_{F_A}(\lambda)}
\text{ and}\\
\EE|X\cap A| &= \sum_{A_0\subset A} \frac{\lambda Z_{G[A_0]}'(\lambda)}{Z_{G[A_0]}(\lambda)} \Pr(A\setminus N(X\setminus A) = A_0) = \EE\frac{\lambda Z_{F_A}'(\lambda)}{Z_{F_A}(\lambda)}, 
\end{align*}
as required for the corollary.
\medskip

We note for $u\in V$ that, by taking $A=\{u\}$, we can also deduce from the spatial Markov property that
\begin{align*}
\Pr(u\in X \mid X\cap N(u)=\varnothing) = \frac{\lambda}{1+\lambda},
\end{align*}
which is as claimed in Section~\ref{sec:intro}. It follows from~\eqref{eqn:hardcore} with $A=N(u)$ that
\begin{align*}
\Pr(u\in X) = \frac{\lambda}{1+\lambda}\EE\frac{1}{Z_{F_{N(u)}}(\lambda)}.
\end{align*}

Therefore, what equation~\eqref{eqn:hardcore} suggests is that in order to establish conditions of the form in~\eqref{eq:lococcgreedy}, it is useful to understand the {\em local} behaviour of the hard-core model. In particular, we want to understand $Z_F(\lambda)$ for induced neighbourhood subgraphs $F$ in $G$, that is, for the induced subgraphs $F$ of $G[N(u)]$ for $u\in V$.
These remarks motivate the following definition, which is central to the chapter.

\begin{definition}
Given $\lambda,\beta,\gamma>0$, if $X$ has the hard-core distribution on $G$ at fugacity $\lambda$, then we say
 $X$ satisfies {\em local $(\beta,\gamma)$-occupancy} if the following holds. For each $u\in V$,
\begin{align}\label{eq:lococc}
\beta \frac{\lambda}{1+\lambda}\frac{1}{Z_{F}(\lambda)} + \gamma\frac{\lambda Z_{F}'(\lambda)}{Z_{F}(\lambda)} \ge 1
\end{align}
for each induced subgraph $F$ of the subgraph $G[N(u)]$ induced by the neighbourhood $N(u)$ of $u$.
\end{definition}


In the next section, we give some examples of how  we can find more tailored choices for $(\beta,\gamma)$ than the ones in Section~\ref{sec:intro}, given various local structural assumptions on $G$.

\begin{lemma}\label{lem:lococc}
If $X$ has local $(\beta,\gamma)$-occupancy, 
then the following hold.
\begin{enumerate}
\item\label{itm:lococcfractional} {\em (Local occupancy for fractional colouring)} For each subset $B\subset V$, there is a distribution over the collection $\cI(G[B])$ of independent sets of the subgraph $G[B]$ of $G$ induced by $B$ such that, writing $X_B$ for the random independent set, it holds for all $u\in B$ that
\begin{align*}
\beta\Pr(u\in X_B) + \gamma\EE|X_B\cap N_{G[B]}(u)|\ge 1. 
\end{align*}
\item\label{itm:lococcoccfrac} {\em (Local occupancy for occupancy fraction)} If $u$ is a vertex chosen independently and uniformly at random from $V$, then
\begin{align*}
\EE_u\left[ \beta\Pr(u\in X) + \gamma\EE|X\cap N(u)| \right] \ge 1. 
\end{align*}
\end{enumerate}
\end{lemma}

For part~\ref{itm:lococcfractional}, we choose the distribution of $X_B$ as the hard-core model on $G[B]$ at fugacity $\lambda$. 
With $F$ chosen as $G[N_{G[B]}(u)]$, it is an induced subgraph of $G[N(u)]$. Since~\eqref{eqn:hardcore} and the accompanying remarks also apply to $X_B$ by our assumption on $X$, we derive the inequality immediately from condition~\eqref{eq:lococc}.
Part~\ref{itm:lococcoccfrac} follows similarly from~\eqref{eqn:hardcore} and the accompanying remarks, where we take $F$ to be the entire neighbourhood subgraph $G[N(u)]$, and then average over the choice of $u$.
\medskip

Of the conditions of the Key Definition, 
Theorem~\ref{lem:lococc}\ref{itm:lococcfractional} and Theorem~\ref{lem:lococc}\ref{itm:lococcoccfrac}, the first is the strongest while the third is in a sense the weakest.
However, in practice our local optimization analyses in the hard-core model do not yield any different values for $\beta$ and $\gamma$ when we target these conditions separately.
In fact, it is the resemblance between these conditions that lies at the root of this chapter.
Thus, in an obvious abuse, we may refer to {\em any} of these three conditions as verifying local $(\beta,\gamma)$-occupancy for the distribution of $X$.

The reason, however, for making explicit these three versions of local occupancy is that they lead to three different implications for the global (independent set) structure of $G$. Moreover, whereas we have assumed throughout the section that $X$ has the hard-core distribution, the assertions in parts~\ref{itm:lococcfractional} and~\ref{itm:lococcoccfrac} of Theorem~\ref{lem:lococc} are the essential hypotheses for the conclusions of the following two results. In other words,  for these results, we do not in principle need that $X$ be distributed according to the hard-core model.

Let us now see how the weakest of these three conditions yields a lower bound on $\EE |X|$, the expected number of vertices of $G$ occupied by $X$. We frequently refer to $\EE |X|/|V|$ as the {\em occupancy fraction}, and we write $\alpha_G(\lambda)$ for this value when $X$ has the hard-core distribution. Note also that the independence number $\alpha(G)$ of $G$ satisfies $\alpha(G)\ge \EE |X|$, by the probabilistic method.

\begin{lemma}[Occupancy fraction via local occupancy]\label{lem:occfrac}
If $X$ satisfies local $(\beta,\gamma)$-occupancy for some $\beta,\gamma>0$ in the sense of Theorem~\ref{lem:lococc}\ref{itm:lococcoccfrac}, 
then \[\frac{\EE |X|}{|V|} \ge \frac{1}{\beta+\gamma \Delta}.\]
\end{lemma}

To see this result, note that, by the condition in Theorem~\ref{lem:lococc}\ref{itm:lococcoccfrac}, the sum of
\[
\beta\Pr(u\in X) + \gamma\EE|X\cap N(u)|
\]
over all $u\in V$ is at least $|V|$. Moreover, this sum is
\begin{align*}
\beta\sum_u\Pr(u\in X) + \gamma\sum_u\sum_{v\sim u}\Pr(v\in X)
\le \beta \EE|X| + \gamma \Delta \EE|X|,
\end{align*} 
where the inequality holds because each vertex $v$ appears at most $\Delta$ times in the double sum. Rearranging then yields the desired result.
\medskip

While the occupancy fraction bound is a simple application of local occupancy, it is a crucial aspect of the technique, because in some settings we can prove that the given bound is asymptotically tight. This is not necessarily true for applications to the independence number or to the (fractional) chromatic number. 
Another important point is that, since larger values of $\lambda$s bias to larger independent sets, the expectation $\EE|X|$ is monotone increasing in $\lambda$ for any fixed graph (see~\cite[Prop.~1]{DJPR17average}), but the quantity $1/(\beta+\gamma\Delta)$ with which we bound it might not be. 
One reason that the case of triangle-free graphs is more subtle and challenging than the general case described in Section~\ref{sec:intro} is precisely this behaviour of $\beta+\gamma\Delta$ as $\lambda$ varies.

We next show how the second of these three conditions -- the one that most closely matches~\eqref{eq:lococcgreedy} -- yields an upper bound on the fractional chromatic number of $G$.

Although there are several equivalent definitions of fractional colouring, we prefer a probabilistic one, in terms of uniform occupancy. A {\em fractional $k$-colouring} of $G$ is a distribution over $\cI$ for which it holds that $\Pr(u\in X_f) \ge 1/k$ for each $u\in V$, where  $X_f$ denotes the random independent set. Note that by summing this over all $u\in V$, we have $\EE |X_f| \ge |V|/k$. The {\em fractional chromatic number} $\chi_f(G)$ is the least $k$ for which $G$ admits a fractional $k$-colouring. Note also that one can see fractional colouring as a relaxation of usual proper colouring, with $\chi_f(G)\le \chi(G)$ always, by noting the simple fractional $\chi(G)$-colouring that selects uniformly from one of the $\chi(G)$ independent sets induced by an optimal proper colouring of $G$.

\begin{lemma}[Fractional chromatic number via local occupancy]\label{lem:fractional}
If $X$ satisfies local $(\beta,\gamma)$-occupancy for some $\beta,\gamma>0$ in the sense of Theorem~\ref{lem:lococc}\ref{itm:lococcfractional},
then \[\chi_f(G) \le\beta+\gamma \Delta.\]
\end{lemma}

To see how this result holds, let us write $k$ for the desired bound $\beta+\gamma\Delta$.
We proceed to build up a fractional colouring ${w}$ in several iterations, and we write ${w}(I)$ for the probability mass assigned to the independent set $I$ and ${w}(u)$ for the total probability mass $\sum_{I\in\cI, I \ni u}{w}(I)$ assigned to independent sets containing the vertex $u$. 
Through the iterations, ${w}$ is a {\em partial} fractional $k$-colouring of $G$, in the sense of not yet having satisfied the condition that ${w}(u)\ge 1/k$  for all $u\in V$, while the total weight ${w}(G)=\sum_{I\in \cI} {w}(I)$ over $\cI$ remains at most $1$.
The idea is to accrue probability mass to the independent sets according to $X_{V(H)}$, for the relevant induced subgraph $H$ of $G$, until the mass containing some vertex is $1/k$, whereupon we remove that vertex from $H$ and continue. This idea is made precise in Algorithm~\ref{alg-greedy}.
\begin{algorithm}\caption{The greedy fractional colouring algorithm}\label{alg-greedy}
    \begin{algorithmic}[0]
        \FOR{$I\in \cI$}
            \STATE${w}(I)\gets0$
        \ENDFOR
        \STATE $H\gets G$
        \WHILE{$|V(H)|>0$}
            \STATE $\displaystyle \tau \gets \min\left\{ \min_{u\in V(H)} \frac{1/k-{w}(u)}{\Pr(u\in X_{V(H)})}, \min_{u\in V(H)} \frac{\beta+\gamma\deg(u)}{k}-{w}(G) \right\}$
            \FOR{$I\in\cI(H)$}
                \STATE ${w}(I)\gets{w}(I)+\Pr(X_{V(H)}=I) \tau$
            \ENDFOR
            \STATE $H \gets H - \{u \in V(H) : {w}(u)=1/k\}$
        \ENDWHILE
    \end{algorithmic}
\end{algorithm}

\noindent
As far as we are aware, Algorithm~\ref{alg-greedy} was first given in~\cite[pp.~244--246]{molloy_graph_2002}, in a slightly more specialized context, which we also discuss later in this chapter.

It remains to prove that Algorithm~\ref{alg-greedy} builds the desired fractional $k$-colouring of $G$.
For the analysis, we index the iterations: for $i=0,1,\dots$, let $H_i$, ${w}_i$ and $\tau_i$ denote the corresponding $H$, ${w}$ and $\tau$ in the $i$th iteration of the `while' loop of Algorithm~\ref{alg-greedy}, prior to updating the sequence.
Writing $V_i = V(H_i)$ for brevity,
we note that $V=V_0 \supseteq V_1 \supseteq V_2 \supseteq \cdots$.
We also have ${w}_{i+1}(v) = \sum_{k=0}^i \Pr(v\in X_{V_k})\tau_k$ for any $v\in V_i$, and ${w}_{i+1}(G) = \sum_{k=0}^i \tau_k$.

By the choice of $\tau_i$, if there is some $u\in V_i$ (that is, with ${w}_i(u) < 1/k$), then ${w}_{i+1}(G) \le (\beta+\gamma\deg(u))/k \le 1$. So we need only show that Algorithm~\ref{alg-greedy} terminates, and to do so, it suffices to show that $V_i \setminus V_{i+1} \ne \varnothing$ for all $i$.
Now if
\[\tau_i = \min_{u\in V_i} \frac{1/k-{w}_i(u)}{\Pr(u\in X_{V_i})}\,,
\]
then there must be some vertex $u\in V_i$ for which ${w}_i(u) < 1/k$ and ${w}_{i+1}(u)=1/k$, and so $u\in V_i \setminus V_{i+1}$ and the result follows. We may therefore assume that there is some vertex $u \in V_i$ for which $\tau_i = (\beta+\gamma\deg(u))/k-{w}_i(G)$, and so ${w}_{i+1}(G)=(\beta+\gamma\deg(u))/k$. 

For any $j\in\{0,1,\dots,i\}$, we know by our assumption of the condition in Theorem~\ref{lem:lococc},  part~\ref{itm:lococcfractional}, that
\[ \beta\Pr(u\in X_{V_j}) + \gamma\EE|X_{V_j}\cap N_{H_j}(u)| \geq 1\,,\] and so 
\[ \beta\Pr(u\in X_{V_j})\tau_j + \gamma\sum_{z\in N_{H_j}(u)} \Pr(z\in X_{V_j})\tau_j \ge \tau_j.\]
By summing the last inequality over all such $j$ and using the above identities for ${w}_{i+1}$, we obtain
\begin{align*}
\frac{\beta+\gamma\deg(u)}{k} = {w}_{i+1}(G) \le
 \beta {w}_{i+1}(u) + \gamma \sum_{z\in N_G(u)} {w}_{i+1}(z).
\end{align*}
Note that the choice of $\tau_i$ ensures that the weights ${w}_{i+1}(z)$ do not exceed $1/k$, and so the right-hand side above is at most
\begin{align*}
\beta {w}_{i+1}(u) + \frac{\gamma \deg(u)}{k}.
\end{align*}
This gives $\beta {w}_{i+1}(u) \geq \beta/k$, and so ${w}_{i+1}(u)=1/k$. 
This means that $u\in V_i \setminus V_{i+1}$, as required for Theorem~\ref{lem:fractional}.
\medskip

Later in this chapter, we show how the strongest of the local occupancy conditions, in the Key Definition, 
together with some mild additional technical conditions, yields an upper bound of roughly $\beta+\gamma\Delta$ on the (list) chromatic number $\chi(G)$ of $G$.

In all three cases (the occupancy fraction, the fractional chromatic number and the chromatic number), we reduce the main task of finding good bounds to the following minimization:
\begin{align*}
\text{minimize} & \qquad \beta+\gamma\Delta\\
\text{subject to} & \qquad \text{$X$ satisfies local $(\beta,\gamma)$-occupancy}
\end{align*}
It is natural to carry out such an optimization task under some local sparsity assumption on $G$, and in the next section we give a few representative examples.

\section{Local occupancy examples}\label{sec:lococcanalysis}

We begin this section by pointing out that the strongest form of local occupancy holds in the general setting discussed in Section~\ref{sec:intro}. 

\begin{theorem}\label{thm:generallo}
    For every $\lambda>0$, the hard-core model at fugacity $\lambda$ on any graph $G$ satisfies local $(\beta,\gamma)$-occupancy with $\beta= 1+ 1/\lambda$ and $\gamma=1$. 
\end{theorem}
    To see this, let $u$ be a vertex of $G$ and let $F\subset G[N(u)]$. Then with $\beta = 1+1/\lambda$ and $\gamma=1$, we want to show that $1 + \lambda Z'_F(\lambda) \ge Z_F(\lambda)$.
    But this follows from the fact that $Z_F(\lambda)$ is a polynomial with non-negative coefficients such that $Z_F(0)=1$, establishing the theorem. 
    Note that equality holds only when $Z_F(\lambda)$ is of degree at most 1.
\medskip

By Theorem~\ref{lem:occfrac}, we can also conclude from Theorem~\ref{thm:generallo} that the hard-core model on $G$ has occupancy fraction at least $1/(1+1/\lambda+\Delta)$, where $\Delta$ is the maximum degree of $G$. As we noted in Section~\ref{sec:intro}, this bound is attained for the complete graph $K_{\Delta+1}$, and so no choice of $(\beta,\gamma)$ in Theorem~\ref{thm:generallo} yields a  strictly smaller value of $\beta+\gamma\Delta$.

\subsection*{Triangle-free graphs}

Further examples of local occupancy require the \emph{Lambert $W$-function}, which is the inverse of the function $x\mapsto xe^x$. We exclusively use the real branch of $W$ defined on $[-1,\infty)$ and require the basic property that $W(x)=\log x-\log\log x+o(1)$ as $x\to\infty$ (see, for example,~\cite{corless_lambert_1996}).

\begin{theorem}\label{thm:tflococc}
    For any $\lambda>0$ and $\gamma>0$, the hard-core model at fugacity $\lambda$ on any triangle-free graph satisfies local $(\beta,\gamma)$-occupancy when 
    \[ \beta = \frac{\gamma{(1+\lambda)}^{(1+\lambda)/(\gamma\lambda)}}{e\log(1+\lambda)}. \]
    
    Moreover, for $d>0$, the choice 
    \[ \gamma  = \frac{1+\lambda}{\lambda}\frac{\log(1+\lambda)}{1+W(d\log(1+\lambda))} \]
    minimizes $\beta+\gamma d$ and yields
    \[ \beta+\gamma d = \frac{1+\lambda}{\lambda} \frac{d\log(1+\lambda)}{W(d\log(1+\lambda))}. \]
    In the limit as $d\to\infty$, provided that $\lambda=o(1)$ with $\log(1/\lambda) = o(\log d)$ (such as when $\lambda=1/\log d$), we have $\beta+\gamma d\sim d/\log d$.
\end{theorem}

It is important to note that the above expression for $\beta+\gamma d$ is not monotone in $\lambda$. If we are interested in the best asymptotic value as $d\to\infty$, then we are forced to consider $\lambda\to 0$ so that the factor $\frac{1+\lambda}{\lambda}\log(1+\lambda)$ is $1-o(1)$, but we also need that $\lambda$ does not vanish so fast so that $W(d\log(1+\lambda)) = o(\log d)$.
Any choice of $\lambda=o(1)$ with $\log(1/\lambda) = o(\log d)$ suffices, such as $\lambda = 1/\log d$.

Let us now see why Theorem~\ref{thm:tflococc} holds.
    In a triangle-free graph $G$, the neighbourhood of any vertex is an independent set, and so it suffices to verify the condition of the Key Definition 
    when $F$ is an edgeless graph on $y$ vertices, $0\le y\le \Delta$. 
    In this case, $Z_F(\lambda) = (1+\lambda)^y$ and $\lambda Z_F'(\lambda)/Z_F(\lambda)=y\lambda/(1+\lambda)$. 
    Thus, the left-hand side of the inequality~\eqref{eq:lococc} is the smooth function 
    \[ g(y) = \frac{\lambda}{1+\lambda}\left(\frac{\beta}{(1+\lambda)^y} + \gamma y\right). \]
    It is straightforward to verify that when $\beta\ge0$, $g$ is convex on the domain $y\ge 0$, because 
    \[ g''(y) = \frac{\beta\lambda}{(1+\lambda)^{1+y}}(\log(1+\lambda))^2 \ge 0. \]
    It follows that $g(y)$ is minimized at its unique stationary point $y^*$, which occurs when 
    \[ g'(y^*) = \frac{\lambda}{1+\lambda}\left(\gamma - \frac{\beta\log(1+\lambda)}{(1+\lambda)^{y^*}}\right) =0 \quad\text{ or, equivalently, }\quad y^* = \frac{\log\left(\beta\log(1+\lambda)/\gamma\right)}{\log(1+\lambda)}.\]
    Since $g(y)\ge g(y^*)$, we can solve for $\beta$ in $g(y^*)=1$ to obtain the stated $\beta=\beta(\gamma,\lambda)$.

    With $\beta$ as above, and given any $d>0$, we now use the same method to find $\gamma$ which minimizes $h(\gamma) = \beta+\gamma d$. 
    We observe that $h$ is convex on the domain $\gamma>0$, and we solve for the stationary point. Write $\theta = \frac{1+\lambda}{\gamma\lambda}\log(1+\lambda)-1$, so that $\theta>-1$ and
    \[ \frac{\partial \theta}{\partial\gamma} = - \frac{(1+\lambda)\log(1+\lambda)}{\gamma^2\lambda}. \]
    We now compute
    \begin{align*}
        h'(\gamma) &= d+\frac{(1+\lambda)^{(1+\lambda)/(\gamma\lambda)} (\gamma \lambda-(1+\lambda) \log (1+\lambda))}{e \gamma \lambda \log (1+\lambda)} = d - \frac{\theta e^{\theta}}{\log(1+\lambda)},
        \\h''(\gamma) &= \frac{1+\lambda}{\gamma^2\lambda}(1+\theta)e^\theta  >0. 
    \end{align*} 
    To obtain the desired $\gamma=\gamma(d,\lambda)$,
    we note that the stationary point $h$ occurs when
    \[ d\log(1+\lambda) = \theta e^\theta  \quad\text{ or, equivalently, }\quad  \theta = W(d\log(1+\lambda)).\]
    The final assertion follows from the fact that $W(x) = \log x - \log \log x + o(1)$ as $x\to\infty$.
    This concludes the proof of Theorem~\ref{thm:tflococc}.
\medskip

With this result, we can now deduce the following corollary.

\begin{corollary}
Let $\lambda>0$.
    If $G$ is a triangle-free graph of maximum degree $\Delta$, then
    \begin{align*}
    	\alpha_G(\lambda) \ge (1-o(1))\frac{\log\Delta}{\Delta}
        \quad\text{and}\quad \chi_f(G) \le (1+o(1))\frac{\Delta}{\log\Delta},
    \end{align*}
    where the $o(1)$ terms are as $\Delta\to\infty$.
\end{corollary}
Both of these bounds follow from the asymptotic analysis of $\beta+\gamma\Delta$ in Theorem~\ref{thm:tflococc}. 
    We obtain the first bound from Theorem~\ref{lem:occfrac} and the monotonicity of $\alpha_G(\lambda)$ in terms of $\lambda$,
    and the second bound follows from Theorem~\ref{lem:fractional}. 
\medskip

\subsection*{Locally sparse graphs}

The main reason that we are interested in triangle-free graphs is that Theorem~\ref{thm:generallo} is tight in disjoint unions of complete graphs, and triangle-freeness is a natural and simple condition that forces a graph to be maximally `far' from a disjoint union of complete graphs, at least locally. 
In the above analysis, the lack of triangles as subgraphs manifests itself most clearly as the condition that the neighbourhoods are independent sets. 
There are many ways to weaken this condition that still result in a significant departure from the behaviour of complete graphs, and we highlight a couple of natural ways for which local occupancy yields good results. 
Both of our settings are efficiently expressed via an upper bound on the average degree of the graphs $F$, which can appear in the condition of the Key Definition; 
the triangle-free case corresponds to enforcing  the average degree of such graphs $F$ to be $0$.

We start with a basic analysis for graphs with bounded average degree.

\begin{lemma}\label{lem:hcmavgdeg}
    For any graph $G$ on $n$ vertices with average degree at most $d$, and any $\lambda>0$,
    \begin{align}
        \frac{\lambda Z_G'(\lambda)}{Z_G(\lambda)}&\ge \frac{\lambda}{1+\lambda}n(1+\lambda)^{-d}\text{ and}\label{eq:EXlbavgdeg}
        \\\log Z_G(\lambda) &\ge
        \begin{cases}
            n\log(1+\lambda) &d=0,\\
            (n/d)\left(1-(1+\lambda)^{-d}\right) & d>0.
        \end{cases}\nonumber
    \end{align}
\end{lemma}
To see this result, let $X$ be a random independent set in $G$ from the hard-core model at fugacity $\lambda$. For any $u\in V(G)$, the spatial Markov property gives 
\[ \Pr(u\in X) = \frac{\lambda}{1+\lambda}\Pr(N(u)\cap X=\varnothing) = \frac{\lambda}{1+\lambda}\EE[1/Z_{F_{N(u)}}(\lambda)]\ge \frac{\lambda}{1+\lambda}(1+\lambda)^{-d_u},\]
where $F_{N(u)} = G[N(u)\setminus(X\setminus N(u))]$ and $d_u$ is the degree of $u$ in $G$. 
The final inequality holds because any realization of $F_{N(u)}$ has $Z_{F_{N(u)}}(\lambda) \le (1+\lambda)^{d_u}$.
We obtain the inequality~\eqref{eq:EXlbavgdeg} from convexity:
\[ \frac{\lambda Z_G'(\lambda)}{Z_G(\lambda)} = \EE|X| \ge \frac{\lambda}{1+\lambda}\sum_{u \in V(G)} (1+\lambda^{-d_u}) \ge \frac{\lambda}{1+\lambda}n(1+\lambda)^{-d}. \]
Recalling~\eqref{eq:Z'}, the second inequality follows from integrating both sides of the first, establishing Theorem~\ref{lem:hcmavgdeg}.
\medskip

Using the above result, we give a generalization of Theorem~\ref{thm:tflococc} whose guarantee degrades smoothly as the average degree bound increases. 
In the following subsections, we apply this result to two natural settings.

\begin{theorem}\label{thm:lomad}
    For any $\avdeg>0$ and $\lambda>0$, consider the hard-core model at fugacity $\lambda$ on any graph $G$ for which the average degree of any subgraph induced by a subset of a neighbourhood in $G$ is at most $\avdeg$.  Then it satisfies, for any $\gamma>0$, local $(\beta,\gamma)$-occupancy when 
    \[ \beta = \frac{\gamma{(1+\lambda)}^{(1+\lambda)^{1+\avdeg}/(\gamma\lambda)-\avdeg}}{e\log(1+\lambda)}. \]
    
    Moreover, for $d>0$, the choice 
    \[ \gamma  = \frac{1+\lambda}{\lambda}\cdot\frac{(1+\lambda)^{\avdeg}\log(1+\lambda)}{1+W(d(1+\lambda)^{\avdeg}\log(1+\lambda))} \]
    minimizes $\beta+\gamma d$ and yields
    \[ \beta+\gamma d = \frac{1+\lambda}{\lambda} \frac{d(1+\lambda)^{\avdeg}\log(1+\lambda)}{W(d(1+\lambda)^{\avdeg}\log(1+\lambda))}. \]
    In the limit as $d\to\infty$, provided that $\lambda=o(\avdeg)$ with $\log(\avdeg/\lambda) = o(\log d)$ (such as when $\lambda=1/(\avdeg\log d)$), we have $\beta+\gamma d\sim d/\log(d/\avdeg)$.
\end{theorem}
    To prove this, it suffices to verify the condition of the Key Definition 
    when $F$ is a graph on $y$ vertices, $0\le y\le \Delta$, with average degree at most $\avdeg$. Then, by Theorem~\ref{lem:hcmavgdeg}, $Z_F(\lambda) \le (1+\lambda)^y$ and $\lambda Z_F'(\lambda)/Z_F(\lambda)\ge y\lambda/(1+\lambda)^{1+\avdeg}$. 
    Thus, the left-hand side of the inequality~\eqref{eq:lococc} is at least $g(y)$, the smooth function given by
    \[ g(y) = \frac{\lambda}{1+\lambda}\left(\frac{\beta}{(1+\lambda)^y} + \gamma y (1+\lambda)^{\avdeg}\right). \]
    We proceed as in the triangle-free case, proving that when $\beta,\avdeg>0$, $g$ is convex on the domain $y\ge 0$. 
    We then solve for the stationary point $y^*$, and on setting $g(y^*)=1$ we solve for $\beta$ to obtain the stated function. 
    As before, given any $d>0$, we can find $\gamma$ which minimizes $h(\gamma) = \beta+\gamma d$ with the same method. 
    When $\avdeg\ge 0$, $h$ is convex on the domain $\gamma>0$ and the given $\gamma$ is the stationary point of $h$.
    The asymptotic properties of the function $W$ again yield the final statement, establishing Theorem~\ref{thm:lomad}.
\medskip

\subsection*{Bounded local triangle count}

A natural weakening of the triangle-free condition is that the graph contains few triangles, and we might expect that the bounds we obtained for triangle-free graphs degrade smoothly as the number of triangles increases; this has already been pursued with respect to the independence number in~\cite{AKS81} and~\cite{She83}. 

Given the way that local occupancy requires only an understanding of neighbourhoods, it is natural to consider a constraint $t$ on the number of triangles containing any given vertex $u$. 
It turns out that the average degree of any graph $F$ arising as a subgraph of a neighbourhood $G[N(u)]$ can be bounded in terms of $t$, as follows. 

\begin{lemma}\label{lem:densitylocallysparse}
    For any graph $G$ in which each vertex is contained in at most $t$ triangles and any $u\in V(G)$, 
     the average degree of any subgraph $F\subset G[N(u)]$ is at most $\sqrt{2t}$.
\end{lemma}
To see this, note that, if $F$ has $y$ vertices, then its average degree is at most $y-1$. 
Similarly, $F$ has at most $t$ edges since each edge corresponds to a triangle in $G$ containing $u$, and so the average degree of $F$ is at most $2t/y$.
Balancing these bounds gives an average degree of at most $(\sqrt{8t+1}-1)/2 \le \sqrt{2t}$, as required for Theorem~\ref{lem:densitylocallysparse}.
\medskip

With this result, we can deduce the following corollary from Theorems~\ref{lem:occfrac} and~\ref{lem:fractional}.

\begin{corollary}
Let $\lambda>0$.
    Let $G$ be a graph of maximum degree $\Delta$ such that each vertex is contained in at most $t$ triangles, for some $t\ge 0$ satisfying $ t= o(\Delta^2)$. Then     
    \begin{align*}
    	\alpha_G(\lambda) \ge (1-o(1))\frac{\log(\Delta/\sqrt{1+t})}{\Delta}
        \quad\text{and}\quad  \chi_f(G) \le (1+o(1))\frac{\Delta}{\log(\Delta/\sqrt{1+t})},
    \end{align*}
    where the $o(1)$ terms are as $\Delta\to\infty$.
\end{corollary}
    Here, both of these bounds follow from the results of Section~\ref{sec:balancing}, given Theorem~\ref{thm:lomad}, the above average degree bound (Theorem~\ref{lem:densitylocallysparse}), and some asymptotic analysis.
    With $1/\lambda = \sqrt{1+t}\cdot \log(\Delta/\sqrt{1+t})$, we have $\lambda=o(1)$, $\log(\sqrt{1+t}/\lambda) = o(\log(\Delta/\sqrt{1+t}))$, and $(1+\lambda)^{\sqrt{2t}} = 1+o(1)$ as $\Delta\to\infty$. With $d=\Delta$, we take $\beta$ and $\gamma$ from Theorem~\ref{thm:lomad}, and obtain 
    \[ \beta+\gamma\Delta \le (1+o(1)) \frac{\Delta}{W((1+o(1))\Delta\lambda)} \le (1+o(1))\frac{\Delta}{\log(\Delta/\sqrt{1+t})}. \]
\medskip

With more fussiness, we can derive results in the range $t=\Theta(\Delta^2)$, although as $t$ approaches $\binom{\Delta}{2}$, alternative methods provide a better understanding; see, for example,~\cite{molloy_graph_2002} and~\cite{HDK22}.

An independence number bound of the form in the above corollary 
was used in~\cite{krivelevich_lower_2004} to establish a bound on the sphere packing density of the form $\theta(d)\ge \Omega(d2^{-d})$. 
An analogue of these ideas in the continuous setting is the main advance of~\cite{jenssen_hard_2019}, resulting in a stronger lower bound on $\theta(d)$ in terms of the leading constant, and deeper problem-specific insights.

\subsection*{Excluded cycle length}

Another natural way to generalize the triangle-free condition is to consider $C_k$-free graphs, with $k\ge 3$, of maximum degree $\Delta$  (where $C_k$ denotes a cycle on $k$ vertices). 
It is not especially clear how bounds on the chromatic number in, say, $C_{8}$-free graphs should be related to those in triangle-free graphs, but nonetheless  from local occupancy we get similar bounds in the limit as $\Delta\to\infty$.

It is impossible to fully express $C_k$-freeness with a condition that applies exclusively in neighbourhoods, but it is certainly true that in $C_k$-free graphs the longest path in a neighbourhood can be on at most $k-2$ vertices. 
This is because if there is a path $P$ on $\ell$ vertices in $G[N(u)]$, then together with $u$ we find a cycle of length $\ell+1$ in $G$ containing $u$. 
It turns out that this property alone permits good local occupancy parameters.

\begin{lemma}\label{lem:densityCkfree}
    Let $k\ge 3$. For any $C_k$-free graph $G$ and any vertex $u\in V(G)$, 
    the average degree of any subgraph $F\subset G[N(u)]$ is at most $k-3$.
\end{lemma}
To see this, note that the longest path in $F$ has at most $k-3$ edges since $G$ is $C_k$-free. 
A theorem of Erd\H{o}s and Gallai~\cite{erdos_maximal_1959} states that the average degree of $F$ is then at most $k-3$, establishing Theorem~\ref{lem:densityCkfree}.
\medskip

With this result, we can deduce the following corollary from Theorems~\ref{lem:occfrac} and~\ref{lem:fractional}.

\begin{corollary}
Let $\lambda>0$.
    Let $G$ be a graph of maximum degree $\Delta$ which is $C_k$-free for some $k\ge 3$ satisfying $k= o(\Delta)$ as $\Delta\to\infty$. Then
    \begin{align*}
    	\alpha_G(\lambda) \ge (1-o(1))\frac{\log(\Delta/k)}{\Delta}
        \quad\text{and}\quad \chi_f(G) \le (1+o(1))\frac{\Delta}{\log(\Delta/k)},
    \end{align*}
    where the $o(1)$ terms are as $\Delta\to\infty$.
\end{corollary}
    Here, both of these bounds follow from the results of Section~\ref{sec:balancing} given Theorem~\ref{thm:lomad}, the above average degree bound (Theorem~\ref{lem:densityCkfree}), and some asymptotic analysis. Set $1/\lambda = k \log(\Delta/k)$ so that, as $\Delta\to\infty$, we have 
    $\lambda=o(1)$, $\log(k/\lambda) = o(\log(\Delta/k))$, and $(1+\lambda)^{k-3}=1+o(1)$. 
    With $d=\Delta$, we take $\beta$ and $\gamma$ from Theorem~\ref{thm:lomad}, and obtain
    \[ \beta+\gamma\Delta \le (1+o(1))\frac{\Delta}{W((1+o(1))\Delta\lambda)} \le (1+o(1))\frac{\Delta}{\log(\Delta/k)}. \] 
\medskip

A more detailed analysis of this method gives a bound for up to $k<\Delta$, but in this range the bounds degrade quickly.
As a way of expressing $C_k$-freeness, the fidelity of the condition that neighbourhoods are $P_{k-1}$-free (where $P_{k-1}$ denotes a path on $k-1$ vertices) would appear to diminish as $k$ grows. In particular, we do not expect these methods to be particularly sharp when $k=\Omega(\Delta)$. 
It would be interesting to extend local occupancy in ways that allow us to capture $C_\Delta$-freeness, although such methods would need to be somewhat less local. 

\subsection*{Fractional version of Reed's conjecture}

A longstanding and difficult conjecture of Reed~\cite{Ree98} posits that for any graph $G$ of clique number $\omega$ and maximum degree $\Delta$, the chromatic number $\chi$ satisfies $\chi\le \lceil(\omega+\Delta+1)/2\rceil$. 
Elegant partial progress on this conjecture was established by Molloy and Reed~\cite{molloy_graph_2002}, who proved that this bound holds for $\chi_f$ without the rounding up.

Their method can be expressed in terms of local occupancy. They showed that if $X$ is a random \emph{maximum} independent set in a graph $G$ of maximum degree $\Delta$ and clique number $\omega$, then for all vertices $u$ we have
\begin{equation}\label{eq:reedlo} 
\tfrac{1}{2}(\omega+1)\Pr(u\in X) + \tfrac{1}{2}\EE|X\cap N(u)| \ge 1.
\end{equation}
As $\lambda\to\infty$, the hard-core model approaches the uniform distribution on maximum independent sets, so this is a local $(\frac12(\omega+1),\frac12)$-occupancy condition of the form in Theorem~\ref{lem:lococc}\ref{itm:lococcfractional} in the case $\lambda=\infty$.
This setting is interesting because the stronger condition in the Key Definition 
does not hold with these parameters. 
Indeed, in the unavoidable case that a subgraph $F\subset G[N(u)]$ has independence number $1$, we have $\Pr(u\in X)=0$ and $\EE|X\cap N(u)|=1$, which makes the left-hand side of inequality~\eqref{eq:lococc} equal to $\frac12$.
This precludes the condition in the Key Definition, 
but if the probability that $\alpha(F_{N(u)})=1$ is sufficiently small, then we can still hope for condition~\eqref{eq:reedlo}.

Molloy and Reed's proof proceeds by revealing $F_{N[u]} = G[N[u]\setminus (X\setminus N[u])]$, a random subgraph of $G$ that always contains $u$ itself, in contrast with the methods discussed so far.  (We have written $N[u]$ to signify the closed neighbourhood $N(u)\cup \{v\}$ of $u$.)

We demonstrate another proof of the inequality~\eqref{eq:reedlo} to highlight a difference between the conditions in the Key Definition 
and Theorem~\ref{lem:lococc}\ref{itm:lococcfractional} and some flexibility in the local occupancy methods.

As in the case of locally sparse graphs, we reveal $F_{N(u)} = G[N(u)\setminus (X\setminus N(u))]$, and then compute the relevant probability and expectation in each case. In the limit $\lambda\to\infty$, the computations are rather simple. The vertex $u$ is occupied if and only if $F_{N(u)}$ is empty, and we always have $|X\cap N(u)| = \alpha(F_{N(u)})$. 
Note that $F_{N(u)}$ is empty if and only if $\alpha(F_{N(u)})=0$, and so we can describe the required quantities in terms of the random variable $A= \alpha(F_{N(u)})$.

Let $p_i = \Pr(A=i)$. In terms of $A$, we have 
\begin{align*} 
\tfrac{1}{2}(\omega+1)&\Pr(u\in X) + \tfrac{1}{2}\EE|X\cap N(u)| 
= \tfrac{1}{2}(\omega+1)p_0 + \tfrac{1}{2}\sum_{i\ge 1}ip_i
\\&\ge \tfrac{1}{2}(\omega+1)p_0 + \tfrac{1}{2}(p_1 + 2(1-p_0-p_1))
 = 1 + \tfrac{1}{2}((\omega-1)p_0-p_1). \end{align*}
This inequality follows from Markov's inequality, but we point out that it is also simply a consequence of the linear constraints that must be satisfied by probabilities: $p_i\ge 0$ for all $i$ and $\sum_i p_i=1$.
The bound~\eqref{eq:reedlo} follows provided that $(\omega-1)p_0\ge p_1$, which we now derive from the spatial Markov property. 

The event $A=1$ corresponds to the realizations of $X$ for which $\alpha(F_{N(u)})=1$. 
Since we are working in the case $\lambda=\infty$, when $A=1$ we must have $|X\cap N(u)|=1$, because $X\cap N(u)$ is distributed according to a maximum independent set in $F_{N(u)}$. 
For any such realization of $X$, we can remove the unique occupied neighbour of $u$ and then insert $u$ to obtain another maximum independent set $X'$ in which $u$ is occupied. So under $X'$ instead, we observe that $A=0$.
Each such $X'$ arises from at most $|F_{N(u)}|$ different realizations of $X$, and when $A=1$ it must be that $|F_{N(u)}|\le \omega-1$, since $F_{N(u)}$ must be a clique in the neighbourhood of $u$ if $A=1$. 
Thus, $(\omega-1)p_0 \ge p_1$.

We emphasize that we have expressed a lower bound on a linear combination of $\Pr(u\in X)$ and $\EE|X\cap N(u)|$ as a linear program whose variables describe the distribution of a random variable. 
We exploited graph structure and the spatial Markov property to express linear constraints on the variables, and we then solved the linear program. 
This is the key idea behind the methods of~\cite{DJPR17independent}, which were developed to give upper bounds on the occupancy fraction of the hard-core and monomer--dimer models. 
In particular, the above argument is inspired by the proof of the upcoming Theorem~\ref{thm:kddocc}.

The latest progress in the original {\em integral} version of Reed's conjecture is found in~\cite{HDK22}.

\section{The local lemma and graph colouring}

A proper $q$-colouring of a graph $G=(V_G,E_G)$ is a function $f:V_G\to[q]$ such that $f(u)\ne f(v)$ for each edge $uv\in E_G$. This corresponds to a partition of $V_G$ into independent sets -- namely, the \emph{colour classes} $f^{-1}(\{i\})$ for $i\in[q]$.

While some connections between independent sets and proper colourings are clear, we begin this section with a construction which shows that the proper $q$-colourings of a graph $G=(V_G,E_G)$ correspond precisely to independent sets of size $|V_G|$ in an auxiliary object known as a \emph{cover}. 
This perspective offers additional generality almost for free, and also allows us to capture the notion of \emph{list colouring}. 
Even if we are only interested in traditional vertex-colouring, after colouring a few vertices we need to keep track of lists of available colours, and so the generality of list colouring arises naturally. 

Formally, to define list colouring we consider a \emph{list assignment}; this is a function $\hat L : V_G\to 2^{\mathbb{N}}$ that assigns to each vertex $v$ of $G$ a subset of the natural numbers, which we call the \emph{list} of $v$.
We say that $G$ is (properly) $\hat L$-colourable if there exists a proper colouring $f$ of $G$ such that $f(v)\in \hat L(v)$ for each vertex $v$, and we say that $G$ is \emph{$q$-choosable} if  $G$ is $\hat L$-colourable for all list assignments $\hat L$ with $|\hat L(v)|=q$.
The \emph{list chromatic number} $\chi_\ell(G)$ of $G$ is the least number $q$ for which $G$ is $q$-choosable.
Note that when $f$ is a list colouring, not only does the fibre $f^{-1}(\{i\})$ induce an independent set, but it also contains only vertices that have $i$ in their list.

Given a graph $G=(V_G,E_G)$ and a list assignment $\hat L$, we construct the cover $\mathscr H=(L,H)$ of $G$ as follows.
Formally, $H=(V_H,E_H)$ is a graph on $\sum_{u\in V_G}|\hat L(u)|$ vertices and $L$ is a map $V_G\to V_H$ which partitions $V_H$ into parts that are indexed by vertices of $G$.
The vertex-set $V_H$ is formed as the (disjoint) union of sets $L(u) = \{ (u,i) : i \in \hat L(u)\}$, for each vertex $u\in V_G$. 
The edges of the cover are given by putting a clique on each $L(u)$, and then, for each edge $uv\in E_G$, connecting $(u,i)$ and $(v,i)$ for each $i\in \hat L(u)\cap \hat L(v)$.

A vertex $x=(u,i)$ in $H$ is incident to two kinds of edges: edges to other vertices in the clique on $L(u)$, and edges to vertices in other parts $L(v)$ for $v\in N_G(u)$. 
For convenience, we define $\deg^*(x)$ to be the number of edges incident to $x$ of the latter type, so $\deg^*(x)=|N_H(x)\setminus L(u)|$.  
We think of $\deg^*(x)$ as the \emph{colour degree}, which represents the number of neighbours of $u$ for which the colour $i$ is also available. When there is ambiguity, we sometimes add a subscript to indicate in which cover the colour degree is taken.  

A proper $\hat L$-colouring $f$ of $G$ corresponds to an independent set of size $|V_G|$ in $H$ -- namely, the set $I = \{ (u,f(u)) : u\in V_G\}$. 
Because of the partition of $|V_H|$ into $|V_G|$ cliques given by $L$, the largest possible size of an independent set in $H$ is $|V_G|$, and so we have reduced the problem of that of showing that $G$ is $\hat L$-colourable to showing that $\alpha(H)=|V_G|$.
For convenience, we define an independent set $I$ in the cover $H$ to be \emph{saturating} if it saturates the bound $|I|\le|V_G|$. (Such a set $I$ is also referred to as an \emph{independent transversal}.)
One might hope that saturating independent sets can be found with local occupancy via Theorem~\ref{lem:occfrac}, but this turns out to be a na\"ive hope in general.

\subsection*{Saturating independent sets}

Given the framing of graph colouring in terms of covers, it is natural to wonder about general conditions under which we can find saturating independent sets. 
Intuitively, when colour-degrees are low and the number of available colours is large, we should hope to succeed. 

\begin{lemma}[`The finishing blow']\label{lem:finishingblow}
    Let $G=(V_G,E_G)$ be a graph and let $\mathscr{H}=(L,H)$ be a cover of $G$. 
    Suppose, for some $d\ge 0$, that $\deg^*(x) \le d$ for all $x\in V_H$ and $|L(u)|\ge 2ed$ for all $u\in V_G$. Then $H$ contains a saturating independent set.
\end{lemma}

For this, we use a symmetric form of the Lov\'asz local lemma (see~\cite{ErLo75} and~\cite{She85}).

\begin{LLL}
	Consider a set $\mathcal{E}$ of (bad) events such that, for each $A\in \mathcal{E}$,
	\begin{enumerate}
		\item $\Pr(A) \le p < 1$ and
		\item $A$ is mutually independent of a set of all but at most $d$ of the other events.
	\end{enumerate}
	If $epd\le1$, then with positive probability none of the events in $\mathcal{E}$ occur.
\end{LLL}

    To see how Theorem~\ref{lem:finishingblow} holds,
    let $S$ be formed by selecting, for each $u\in V_G$, a vertex $x\in L(u)$ uniformly at random. 
    For an edge $xy$ in $H$, let $B_{xy}$ be the bad event that $xy\in H[S]$. 
    Note that $S$ is a saturating independent set if and only if no bad event occurs. Also, if $x\in L(u)$ and $y\in L(v)$, then 
    \[ \Pr(B_{xy}) = \frac{1}{|L(u)||L(v)|}, \] 
    and $B_{xy}$ is independent of $B_{zw}$ unless at least one of $z$ and $w$ lies in $L(u)\cup L(v)$. 
    Thus, in the dependency graph, $B_{xy}$ has degree at most $d|L(u)| + d|L(v)|-1$ (excluding $B_{xy}$ itself). 
    By the Lov\'asz local lemma, the proof is complete if, for all $uv\in E_G$ and $xy\in E_H$ with $x\in L(u)$ and $y\in L(v)$, we have  
    \[ ed \frac{|L(u)|+|L(v)|}{|L(u)||L(v)|} \le 1,\]
    but this follows from the lower bound on the list sizes, establishing the result.
\medskip

A version of this result with a worse leading constant was given by Alon~\cite{Alo88}. Using  topological methods, Haxell~\cite{Hax01} proved it with $2e$ replaced by $2$. In fact, one can prove it with $2e$ replaced by $1+o_d(1)$ via nibble methods (see~\cite{RS02}).
For our purposes, the difference between these results is negligible, and so we opt for a statement that has a relatively simple proof in order to keep the arguments self-contained. 
Awkwardly, as shown in~\cite{BH02}, we cannot replace $2ed$ with $d+1$.

While Theorem~\ref{lem:finishingblow} is a good start, it is not so powerful on its own.
For example, we cannot straightforwardly recover the greedy bound, which states that $\chi(G)\le \Delta+1$. 
If we consider a cover for proper $(\Delta+1)$-colourings of the complete graph $K_{\Delta+1}$, then every colour-degree is $\Delta$ and there is no version of Theorem~\ref{lem:finishingblow} for these parameters. 

A more fruitful approach is to show that there is some independent set $I\in \cI(H)$ for which we can apply Theorem~\ref{lem:finishingblow} in the `leftover' we have after removing $I$ and all its neighbours from $H$.
Given a graph $G=(V_G,E_G)$, a cover $\mathscr{H}=(L,H)$ of $G$, and an independent set $I\in \cI(H)$, we define the \emph{leftover} to be the graph $G_I = G - \{u\in V_G : |I\cap L(u)|=1\}$, together with a cover $\mathscr{H}_I=(L_I,H_I)$ of $G_I$, where
 $H_I=H-N[I]$ and $L_I$ is the partitioning map given by $L_I(u) = L(u)\setminus N[I]$.
 
 As a side remark on nomenclature, given the standard use of `uncovered' in the setting of independent sets in $G$, the term \emph{uncovered cover} is perhaps more evocative than `leftover', but we have chosen to sidestep any confusion that this term might cause.

An important strategy for bounding the chromatic number of $G$ is to show that there is some $I\in \cI(H)$ for which the conditions of Theorem~\ref{lem:finishingblow} apply to $\mathscr{H}_I$.

For a locally sparse graph $G$, this strategy has been useful since the 1990s: both Kim~\cite{Kim95a} (for $G$ of girth $5$) and Johansson~\cite{johansson_asymptotic_1996} (for $G$ triangle-free) found an appropriate set $I$ by nibbling, an approach covered in Chapter~\tomchapter.
It was a breakthrough in 2019 when Molloy~\cite{molloy_list_2019} showed that a suitable set $I$ can be found by using entropy compression to analyze an algorithm for generating a random independent set. 
Bernshteyn~\cite{bernshteyn_johanssonmolloy_2019} subsequently showed that Molloy's argument can be cast in terms of analyzing a uniformly random independent set in the cover using the so-called `lopsided' Lov\'asz local lemma. 
Later work~\cite{davies_graph_2020} involving the authors showed that the conditions under which Bernshteyn's analysis applies follow from local occupancy. 

\subsection*{Random independent sets in cover graphs}

The idea of considering a random independent set $Y$ in $H$, distributed according to the hard-core model at fugacity $\lambda$, is the main conceptual breakthrough that connects local occupancy with proper colourings. 
Only for very large $\lambda$ might we expect to gain information about saturating independent sets in $H$, but for the best results we are forced to handle smaller values of $\lambda$ and find a smaller independent set which allows for an application of Theorem~\ref{lem:finishingblow}.
This is the reason for some `mild additional technical conditions' in the colouring aspects of local occupancy mentioned in Section~\ref{sec:intro}.
We make the following standing assumptions for the remainder of this section:

\begin{itemize}
\item $G=(V_G,E_G)$ is a graph of maximum degree $\Delta$;
\item $X$ is a random independent set of $G$ chosen according to the hard-core model at fugacity $\lambda$, and
$X$ has local $(\beta,\gamma)$-occupancy as per the Key Definition; 
\item $\mathscr{H}=(L,H)$ is a cover of $G$ with $|L(v)|=q$ for all $v\in V_G$, and $u\in V_G$ is an arbitrary fixed vertex; and
\item $Y$ is a random independent set of $H$ chosen according to the hard-core model at fugacity $\lambda$, and $Y'$ is $Y$ conditioned on $u\in V(G_Y)$. 
\end{itemize}

This section is focused on showing that, with suitable parameter choices, there is positive probability that the leftover -- that is, the cover $\mathscr H_Y$ together with the graph $G_Y$ -- satisfies the conditions of Theorem~\ref{lem:finishingblow}.
Since we only need to know about $L_Y(u)$ and $\deg^*(x)$ for $x\in L(u)$ when $u\in V(G_Y)$, we can afford to condition on $u\in V(G_Y)$ in the probabilistic analysis. 
That is why we can fix $u$ and consider the random independent set $Y'$, which by the spatial Markov property (Theorem~\ref{lem:spatialMarkov}) is distributed according to the hard-core model on $H'=H-L(u)$.

We start with the observation that local occupancy for $X$ gives us useful information akin to local occupancy for $Y'$ from the perspective of each $x\in L(u)$.

\begin{lemma}[Local occupancy in the cover graph]\label{lem:coverocc}
    Let $x\in L(u)$ and let $F=H'[N_H(x)\setminus N(Y'\setminus N_H(x))]$, where  $H'=H-L(u)$.
    Then 
    \begin{align}\label{eq:pxF}
        \Pr(x\in L_{Y'}(u) \mid F) &= \frac{1}{Z_F(\lambda)}.
    \end{align}
    Moreover, conditioned on any realization of $Y'\setminus N_H(x)$,
    \begin{equation}\label{eq:coverlococc}
        \beta\frac{\lambda}{1+\lambda} \Pr(x\in L_{Y'}(u)) + \gamma\EE|N_H(x)\cap Y'|  \ge 1.
    \end{equation} 
\end{lemma}
    In some sense, the graph $F$ here is the neighbourhood of $x$ uncovered by $Y'$,
    although one must take care with the definition since we condition on $u\in G_Y$ and we do not wish to include vertices of $L(u)$ in this neighbourhood.

    To see how Theorem~\ref{lem:coverocc} holds, we note that $x\in L_{Y'}(u)$ if and only if $Y'\cap N_{H}(x)=\varnothing$, and also that the set $N_H(x)\cap Y'$ of occupied neighbours of $x$ is simply  the set of occupied vertices in $F$. 
    Then the spatial Markov property (Theorem~\ref{lem:spatialMarkov}) gives  
    \begin{align*}
        \Pr(x\in L_{Y'}(u) \mid F) &= \frac{1}{Z_F(\lambda)}\quad\text{and}\quad\EE[|N_H(x)\cap Y'|\mid F]=\frac{\lambda Z_F'(\lambda)}{Z_F(\lambda)},
    \end{align*}
    as required for the equation~\eqref{eq:pxF}. 
    But $H'[N_{H}(x)]$ is isomorphic to an induced subgraph of $G[N(u)]$, and hence $F$ is too, for any realization of $Y'\setminus N_H(x)$. 
    This permits an application of the local occupancy condition that we assume for $X$, and the inequality~\eqref{eq:coverlococc} now follows from the inequality~\eqref{eq:lococc}.
\medskip

    Theorem~\ref{lem:coverocc} clearly requires the additional power of the condition in the Key Definition, 
     in comparison with the weaker conditions in Theorem~\ref{lem:lococc}. 
    Under the weaker conditions for $X$, it is unclear how one could conclude anything useful about the hard-core model on a cover of $G$.
    For the remainder of this section, we make no further reference to $X$, save for the consequences that it has for $Y'$ as per Theorem~\ref{lem:coverocc}.

The following result gives a lower bound on the quantity $\EE|L_{Y'}(u)|$, the expected size of a list in the leftover cover, and implies an upper bound on the quantity $\EE\deg^*_{\mathscr{H}_{Y'}}(x)$, the expected colour degree of any vertex $x$ in the leftover cover.
This shows that on the level of the first moment the setup of Theorem~\ref{lem:finishingblow} is attainable.
We will subsequently establish concentration. 

For clarity, we state the first moment result with no additional conditioning beyond $u\in V(G_Y)$, but crucially we will also observe that they hold conditioned on any realization of $Y'\setminus L(N(u))$. 
Avoiding explicit conditioning in the proof considerably lightens the notation.

\begin{lemma}\label{lem:coverexpect}
Under the standing assumptions,
    \[\EE|L_{Y'}(u)|\ge \frac{1+\lambda}{\beta\lambda}(q-\gamma\Delta)\]
    and, for any $x\in L(u)$ and any $d\ge 0$, conditioned also on $x\in L_{Y'}(u)$, 
    \begin{equation}\label{eq:pcoldegbound}
    \Pr(\deg^*_{\mathscr{H}_{Y'}}(x) \ge d) \le \max_{\substack{F \subset G[N(u)], \\ |V(F)|= d}}\frac{1}{Z_F(\lambda)}.
    \end{equation}
    Moreover, the same bounds hold when additionally conditioning on any realization of $Y'\setminus L(N(u))$.
\end{lemma}

    Given the identity $\EE R = \sum_{d=1}^\infty \Pr(R \ge d)$, which holds for random variables $R$ supported on the non-negative integers with finite mean, it follows immediately from~\eqref{eq:pcoldegbound} that 
    \[ \EE[\deg^*_{\mathscr{H}_{Y'}}(x)] \le \sum_{d=1}^\Delta \max_{\substack{F \subset G[N(u)], \\ |V(F)|= d}}\frac{1}{Z_F(\lambda)}. \]
    Due to the rather general setup, it is challenging to interpret the implications of Theorem~\ref{lem:coverexpect}. 
    We can substitute in bounds on the $Z_F(\lambda)$ terms subject to local information. 
    In the setting discussed in Section~\ref{sec:intro}, we needed to allow for complete graphs $F$ and hence we could only ensure that $Z_F(\lambda)\ge 1+d\lambda$ for $|V(F)|=d$. But in the setting of triangle-free graphs, $F$ contains no edges and we may use the stronger bound $Z_F(\lambda)\ge(1+\lambda)^d$.
    In either case, we see that Theorem~\ref{lem:coverexpect} provides an understanding of the first moment of quantities relevant to Theorem~\ref{lem:finishingblow}.

    To see why Theorem~\ref{lem:coverexpect} holds, the first inequality follows from Theorem~\ref{lem:coverocc}.
    For, if we sum the inequality~\eqref{eq:coverlococc} over all vertices $x\in L(u)$ we get 
    \[ \beta\frac{\lambda}{1+\lambda}\EE|L_{Y'}(u)| + \gamma\sum_{x\in L(u)}\EE|N_H(x)\cap Y'|  \ge q. \]
    Here, the sum over $x\in L(u)$ simply counts the expected number of neighbours $v$ of $u$ (in $G$) which satisfy $|Y'\cap L(v)| = 1$, and this is clearly at most $\Delta$. 
    Rearranging, we have the required bound on $\EE|L_{Y'}(u)|$.

    By Theorem~\ref{lem:coverocc} again, we have the equality~\eqref{eq:pxF}, which states that $\Pr(x\in L_{Y'}(u) \mid F) = 1/Z_F(\lambda)$ when $F=H'[N_H(x)\setminus (Y'\setminus N_H(x))]$. 
    Thus, for the probability bound, it suffices to ensure that whenever $|V(F)|\ge d$, we have that
    $Z_F(\lambda)$ is large enough. 
    But since all such $F$ are isomorphic to subgraphs of $G[N(u)]$, and since $Z_F(\lambda) \le Z_{F'}(\lambda)$ whenever $F\subset F'$, this is guaranteed by the statement of the theorem.
    
    The same bounds also hold subject to the additional conditioning on $Y'\setminus L(N(u))$ because all we need is a property of all subgraphs $F$ that arise as subgraphs of $G[N(u)]$. In particular, Theorem~\ref{lem:coverocc} holds subject to such conditioning, and this establishes Theorem~\ref{lem:coverexpect}.
\medskip

Concentration is a simple matter in our setup: we show that a Chernoff bound applies for the list size in the leftover, and we use a union bound for the colour degrees.

\begin{lemma}\label{lem:coverconc}
\[
    \Pr\left( |L_{Y'}(u)| \le \frac{1+\lambda}{2\beta\lambda}(q-\gamma\Delta)\right) \le \exp\left(-\frac{1+\lambda}{8\beta\lambda}(q-\gamma\Delta)\right).
\]
For any fixed $d$, the probability that there exists some $x\in L_{Y'}(u)$ with $\deg^*_{\mathscr{H}_Y}(x)\ge d$ is at most 
\[ q \max_{\substack{F \subset G[N(u)], \\ |V(F)|= d}}\frac{1}{Z_F(\lambda)}. \]
Moreover, we have the same bounds conditioned on any realization of $Y'\setminus L(N(u))$.
    \end{lemma}
    The first bound is an application of a well-known Chernoff bound for negatively correlated random variables. 
    Since $|L_{Y'}(u)|$ is the sum of the indicator random variables $Z_x$ that $x\in L_{Y'}(u)$, it suffices to show that $Z'_x = 1-Z_x$ over $x\in L(u)$ are negatively correlated. 
    We might expect negative correlation, by the heuristic that if $x=(u,i)$ is in $x\in L_{Y'}(u)$, then $Y'$ cannot signify the colour $i$ being used on $N(u)$, and hence other colours are more likely to be used on $N(u)$ and so not appear in $L_{Y'}(u)$.
    
    It is enough to show, for all $x\in L(u)$ and $S\subset L(u)\setminus\{x\}$, that
    \[
    \Pr(x\notin L_{Y'}(u) \mid S \cap L_{Y'}(u) = \varnothing ) \le \Pr(x\in L_{Y'}(u)),
    \]
    which is equivalent to
    \[
        \Pr\big(\text{$Y'\cap N_{H}(y) \ne \varnothing$ for all $y\in S \mid Y'\cap N_{H}(x)=\varnothing$}\big) \ge \Pr(S\cap L_{Y'}(u)=\varnothing).
    \]
    In the triangle-free case this holds because the sets $N_{H}(x) \setminus L(u)$ and $N_{H}(S)\setminus L(u)$ are disjoint, and in the general case negative correlation can be established with a more involved argument.

    The second bound in Theorem~\ref{lem:coverconc} follows from the inequality~\eqref{eq:pcoldegbound} and a union bound over the $q$ elements $x\in L(u)$.
\medskip

\subsection*{Putting it all together}

\begin{theorem}[Chromatic number via local occupancy]\label{thm:colouring}
    Suppose that the hard-core model on $G$ at fugacity $\lambda$ has local $(\beta,\gamma)$-occupancy, and that for $q>0$ there is some $d\ge \frac{2}{e}\log(2e\Delta^3)$ such that 
    \begin{align}
        \frac{1+\lambda}{\beta\lambda}(q-\gamma\Delta) &\ge 4ed \quad \text{and}\label{eq:expectlb}\\
        \max_{u\in V_G} \max_{\substack{F \subset G[N(u)], \\ |V(F)|= d}}\frac{1}{Z_F(\lambda)} &\le \frac{1}{2eq\Delta^3}.\label{eq:maxcoldegprob}
    \end{align}
    Then $\chi(G)\le q$.
\end{theorem}
    To see this, we use the Lov\'asz local lemma to guarantee an independent set $I$ in $H$ for which the conditions of Theorem~\ref{lem:finishingblow} hold in $\mathscr{H}_I$. 
    For $u\in V_G$, let $A_u$ be the bad event that $u\in V(G_Y)$ and $|L_Y(u)|< 2ed$, and let $B_u$ be the bad event that $u\in V(G_Y)$ and there exists some $x\in L_Y(u)$ with $\deg^*_{\mathscr{H}_Y}(x) > d$.
    We consider the family of events $C_u = A_u\cup B_u$ in the local lemma. Note that $C_u$ is independent of $C_v$ for any $v\in V_G\setminus N^3[u]$. 
    This is because the bad events are all properties of the list $L_{Y}(u)$ of $u$ in the leftover cover, and this is fixed by the intersection of $Y$ with $N_H(L(u))$. If the distance in $G$ between $u$ and $v$ is at least $4$, then $N_H(L(u))$ and $N_H(L(v))$ are disjoint.
    Since $|N^3[u]| \le \Delta^3$, it suffices to show that 
    $\Pr(A_u), \Pr(B_u)\le 1/(2e\Delta^3)$.

    The bound $\Pr(A_u)\le 1/(2e\Delta^3)$ follows from Theorem~\ref{lem:coverconc}, the condition~\eqref{eq:expectlb} and the lower bound on $d$.
    The bound $\Pr(B_u)\le 1/(2e\Delta^3)$ follows from Theorem~\ref{lem:coverconc} and the condition~\eqref{eq:maxcoldegprob}.
    This establishes Theorem~\ref{thm:colouring}.
\medskip

The astute reader may already have noticed that the setup above straightforwardly generalizes to list colouring, and so we can obtain the following results mainly as consequences of the local occupancy analyses of Section~\ref{sec:lococcanalysis}. 
The first of these is Molloy's theorem~\cite{molloy_list_2019}, while the second and third are smooth generalizations of Molloy's theorem as shown in~\cite{davies_graph_2020}.

\begin{corollary}
    Let $G$ be a graph of maximum degree $\Delta$.
    \begin{enumerate}
        \item\label{part:chi:tf} If $G$ is triangle-free, then $\chi_\ell(G)\le (1+o(1))\Delta/\log\Delta$.
        \item\label{part:chi:sparsenbhds} If each vertex of $G$ is contained in at most $t$ triangles for some $t\le o(\Delta^2)$, then $\chi_\ell(G)\le (1+o(1))\Delta/\log(\Delta/\sqrt{1+t})$.
        \item\label{part:chi:cycles} If $G$ is $C_k$-free for $3\le k \le o(\Delta)$, then $\chi_\ell(G)\le (1+o(1))\Delta/\log(\Delta/k)$.
    \end{enumerate}
\end{corollary}
    We established the fractional variants of these upper bounds in Section~\ref{sec:lococcanalysis}. To see the above corollary, what remains is to show that the `mild' additional conditions of Theorem~\ref{thm:colouring} hold, and to study the asymptotic growth of the resulting bounds on the list chromatic number.
    In each case we have a bound $\avdeg$ on the average degree in any subgraph of a neighbourhood in $G$, and we can thus apply Theorems~\ref{lem:hcmavgdeg} and~\ref{thm:lomad}. 
    For convenience, we have stated conditions which guarantee that $\avdeg=o(\Delta)$.

    For part~\ref{part:chi:tf}, the triangle-free case is $\avdeg=0$. 
    Suppose that $F$ is a subgraph of a neighbourhood in $G$ with at least $d$ vertices. Then ${1}/{Z_F(\lambda)} \le {1}/{b^d}$,
    where $b=1+\lambda$. 
    To satisfy inequality~\eqref{eq:maxcoldegprob} we can choose any $d$ such that $d \ge\log_b(2eq\Delta^3)$,
    and since the desired list chromatic number bound follows from a simple greedy argument when $q\ge \Delta+1$, it suffices to set
    $d = \log_b(2e\Delta^4)$.
    It remains to show that there is some $q=(1+o(1))\Delta/\log(\Delta/\avdeg)$ satisfying~\eqref{eq:expectlb}. 
    Rearranging, we want
    \[ q \ge \frac{4ed\lambda}{1+\lambda}(\beta + \Delta'\gamma), \]
    where $\Delta'=(1+\lambda)\Delta/(4ed\lambda)$. Theorem~\ref{thm:tflococc} tells us that if we can choose $\lambda$ such that $\lambda=o(1)$ and $\log(1/\lambda) =o(\log \Delta')$ as $\Delta\to \infty$, then there are suitable choices of $\beta$ and $\gamma$ such that $\beta+\Delta'\gamma \sim \Delta'/\log\Delta'$.
    With the choices $\lambda=1/\log\Delta$ and $d$ as above, we have that $\Delta' = \Theta(\Delta/\log\Delta)$, which gives us what we require. Then 
    \begin{align*}
        \frac{4ed\lambda}{1+\lambda}(\beta + \Delta'\gamma)
        \sim \frac{4ed\lambda}{1+\lambda} \frac{\Delta'}{\log\Delta'}
        = \frac{\Delta}{
        \log\Delta+\log\frac{1+\lambda}{4ed\lambda}
        } \sim \frac{\Delta}{\log\Delta}.
    \end{align*}

    The case $\avdeg>0$ is similar. Without loss of generality, we may assume that $\avdeg\ge 1$. 
    When $F$ is a subgraph of a neighbourhood in $G$ with at least $d$ vertices, Theorem~\ref{lem:hcmavgdeg} implies that ${1}/{Z_F(\lambda)} \le {1}/{b^d}$,
    with $\log b=(1-(1+\lambda)^{-\avdeg})/\avdeg$. 
    Again, we set
    $d = \log_b(2e\Delta^4)$, and minimize the lower bound, \[ \frac{4ed\lambda}{1+\lambda}(\beta + \Delta'\gamma), \]
    that we require on $q$, 
    where $\Delta'=(1+\lambda)\Delta/(4ed\lambda)$. 
    
    Theorem~\ref{thm:lomad} tells us that if we can choose $\lambda$ such that $\lambda=o(1)$ and $\log(\avdeg/\lambda) =o(\log(\Delta'))$ as $\Delta\to \infty$, then there are suitable choices of $\beta$ and $\gamma$ such that $\beta+\Delta'\gamma \sim \Delta'/\log(\Delta'/\avdeg)$.
    Using that $\avdeg\ge 1$ and $\avdeg= o(\Delta)$, with the choices $1/\lambda = \avdeg \log(\Delta/\avdeg)$ and $d$ as above, we deduce that $\Delta' = \Theta(\Delta/\log\Delta)$ and hence that $\log(\avdeg/\lambda) = \log(\log(\Delta/\avdeg)) = o(\log(\Delta'))$, as required. Then 
    \begin{align*}
        \frac{4ed\lambda}{1+\lambda}(\beta + \Delta'\gamma)
        \sim \frac{4ed\lambda}{1+\lambda} \frac{\Delta'}{\log(\Delta'/\avdeg)}
        = \frac{\Delta}{
        \log(\Delta/\avdeg)+\log\frac{1+\lambda}{4ed\lambda}
        } \sim \frac{\Delta}{\log(\Delta/\avdeg)}.
    \end{align*}
    The result now follows because we can take $\avdeg=\sqrt{2t}$ for part~\ref{part:chi:sparsenbhds} and $\avdeg=k-3$ for part~\ref{part:chi:cycles}. The final detail is to argue that $\log(\Delta/\sqrt{2t}) \sim \log(\Delta/\sqrt{1+t})$ and $\log(\Delta/(k-3)) \sim \log(\Delta/k)$ under the present assumptions, and this is simple asymptotic analysis.
\medskip

\section{The occupancy fraction of regular graphs}

In this section, we change tack slightly by seeking upper bounds on the occupancy fraction $\alpha_G(\lambda)$ over families of graphs $G$. 
It is immediately apparent that this is not an interesting problem for bounded-degree graphs, because
\[ \Pr(u \in X)\le \frac{\lambda}{1+\lambda}\Pr(|X\cap N(u)|=0) \le \frac{\lambda}{1+\lambda} \]
 for any vertex $u$ in any graph, and this bound is tight when $u$ is an isolated vertex. 
For this reason, we impose the stronger condition of $d$-regularity instead of merely an upper bound on the degrees. 

As we have already carried out the requisite analysis in the hard-core model, it remains to carefully optimize them in different classes of graphs. 
We start with a result from~\cite{DJPR17independent}, although we give an alternative proof of Perkins that was published in Zhao's survey~\cite{zhao_extremal_2017}.

\begin{theorem}
\label{thm:kddocc}
    For any $\lambda>0$ and any $d$-regular graph $G=(V,E)$, 
    \begin{align*}
	    \alpha_{G}(\lambda)\le \alpha_{K_{d,d}}(\lambda)
	    \quad\text{ and }\quad
	    Z_G(\lambda)^{1/|V(G)|} \le Z_{K_{d,d}}(\lambda)^{1/(2d)}.
    \end{align*}
\end{theorem}
    To see why this result holds, we write $n=|V|$ and recall that $\alpha_G(\lambda) = \frac{1}{n}\sum_{u\in V}\Pr(u\in X)$, where $X$ is distributed according to the hard-core model on $G$ at fugacity $\lambda$.
    We reveal $Y_u=|X\cap N(u)|$ and observe that
    \[ \alpha_G(\lambda) = \frac{1}{n}\sum_{u\in V}\frac{\lambda}{1+\lambda}\Pr(Y_u=0) = \frac{1}{dn}\sum_{u\in V}\EE Y_u, \]
    by a calculation from Section~\ref{sec:intro}. Here, the first equality follows from the spatial Markov property and the second equality uses the fact that each vertex appears $d$ times in the double sum $\sum_{u\in V}\sum_{v\in N(u)}\Pr(v\in X) = \sum_{u\in V}\EE Y_u$.
    We can further simplify our expressions for $\alpha_G(\lambda)$, by defining $Y$ to be $Y_u$ for a uniformly random choice of $u$ and writing $q_i=\Pr(Y=i)$. By linearity of expectation, we have  
    \begin{align*} 
    \alpha_G(\lambda) &= \frac{\lambda}{1+\lambda}\Pr(Y=0) = \frac{1}{d}\EE Y
    =\frac{\lambda}{1+\lambda}q_0 = \frac{1}{d}\sum_{i=1}^d i q_i.
    \end{align*}
    The `objective function' $\alpha_G(\lambda)$ is clearly linear in the graph-dependent variables $q_i$, and we have linear constraints: $q_i\ge 0$, $\sum_i q_i=1$, and the equality above. 

    The main idea now is to relax the above optimization problem, which formally is over all probability vectors $\mathbf q = (q_0,q_1,\dotsc,q_d)$ that arise as the distribution of $Y$ on a $d$-regular graph, to a linear program over all real vectors $\mathbf q \in \mathbb{R}^{d+1}$. 
    For high-quality bounds, we need a definition of $Y$ that reveals enough graph structure to capture the class of interest, and we want to include as many valid linear constraints as possible in order to reduce the difference between the true optimum and the optimum of the relaxation.
    Thus far, we are missing some important constraints. 
    Rather like the constraint used to establish the condition in~\eqref{eq:reedlo}, for $i\ge 1$ we can take any realization of the random independent set $X$ and random vertex $u$ which yields $Y=i$, and form $X'$ by removing any vertex in $X\cap N(u)$. This gives a pair $(X',u)$ for which we would find $Y=i-1$. 
    The important thing here is to take into account both any overcounting and the change in size $|X|=|X'|+1$ which entails a factor of $\lambda$ in the probabilities. 
    Each such $X'$ can be found by removing a vertex from at most $d-i-1$ different $X$s, and each $X$ can have a vertex removed in $i$ ways, yielding the constraint 
    \begin{equation}\label{eq:kddconst} (d-i-1)q_{i-1}\lambda \ge i q_i. \end{equation}
    Intuitively, these constraints help to reduce the maximum of the linear program, as, for $i\ge 2$, making $q_i$ smaller and $q_j$ larger, for some $1\le j<i$, reduces $\alpha_G(\lambda)$ via the expression as $\EE Y/d$.
    
    It remains to solve the relaxed program and find that things work out perfectly: the optimum of the relaxed program is attainable by a graph in our class -- namely, $K_{d,d}$. 
    While standard duality techniques make this somewhat straightforward, for completeness we give an elementary analysis. 
    
    We start by proving that if a vector $\mathbf q$ achieves the maximum, then the inequality~\eqref{eq:kddconst} must hold with equality for $2\le i\le d$;
    this follows by a simple re-weighting argument. 
    If for $i\ge 2$ we have strict inequality, then for a small enough $\eps>0$ we can increase $q_0$ by $\eps$, decrease $q_{i-1}$ by $\left(\frac{d\lambda}{1+\lambda}+i\right)\eps$ and increase $q_i$ by $\left(\frac{d\lambda}{1+\lambda}+i-1\right)\eps$, while maintaining all the constraints. But since we have increased $q_0$ we have also increased the objective value, which is a contradiction. 
    It is now a simple, if tedious, calculation to show that there is precisely one vector $\mathbf q$ that achieves equality in~\eqref{eq:kddconst} for $2\le i\le d$, while satisfying all the other constraints, and that the vector $\mathbf q$ is exactly the distribution of $Y$ when the graph is $K_{d,d}$. 

    Recalling~\eqref{eq:Z'}, we then integrate $\EE|X|/\lambda$ to obtain an estimate on $\log Z_G(\lambda)$ for the second assertion. This establishes Theorem~\ref{thm:kddocc}.
\medskip

The above ideas can be applied to a setting with much more challenging graph structure. 
Recall that, given a graph $G=(V,E)$, the \emph{line graph} $L(G)$ is a graph with vertex-set $E$, where we connect two edges of $G$ in $L(G)$ if they share a common vertex. 
As the line graph of a $d$-regular simple graph is $2(d-1)$-regular, it turns out that if we tolerate a slight change in the parameter we can  then ask for a version of Theorem~\ref{thm:kddocc} over the class of line graphs of regular graphs. 
For $d\ge3$, the graph $K_{d,d}$ is not a line graph, and so we might hope for an improved upper bound. 
The following tight result from~\cite{DJPR17independent} establishes that $L(K_{d,d})$ is the extremal graph.

\begin{theorem}
\label{thm:kddoccmonomer}
    For any $\lambda>0$ and any $d$-regular graph $G=(V,E)$,  
    \begin{align*}
	    \alpha_{L(G)}(\lambda)\le \alpha_{L(K_{d,d})}(\lambda)
	    \quad\text{ and }\quad
	    Z_{L(G)}(\lambda)^{1/|E|} \le Z_{L(K_{d,d})}(\lambda)^{1/d^2}.
    \end{align*}
\end{theorem}

The proof of this result is similar to the above, but in the line graph we take a random independent set $X$ from the hard-core model and reveal the graph structure on $F_{N[e]}$, instead of merely $|X\cap N(e)|$; 
this gives us more room to express constraints on the distributions involved. 
The solution of the resulting linear program is also significantly more involved.
We note that, whereas there are alternative and sometimes more general proofs of the partition function bound in Theorem~\ref{thm:kddocc}, using the entropy method of Kahn~\cite{kahn_entropy_2001},~\cite{zhao_number_2010} or H\"older-type inequalities (see~\cite{lubetzky_replica_2015},~\cite{sah_number_2019} and~\cite{sah_reverse_2020}), we know of no alternative proof of the partition function bound in Theorem~\ref{thm:kddoccmonomer}. 
We can also observe consequences of the additional strength of occupancy fraction bounds over partition function bounds in some circumstances. 
It was shown in~\cite{DJPR17average} that one can integrate the bound on occupancy fraction in triangle-free graphs to get the best-known lower bound on the total number $Z_G(1)$ of independent sets in $n$-vertex triangle-free graphs $G$. 
The analogous result for the hard-sphere model in~\cite{jenssen_hard_2019} gives a lower bound on the entropy (roughly, the quantity) of sphere packings, rather than just their density.

\section{Barriers}\label{sec:barriers}

Our main goal in this chapter has been to show the main principles behind the utility of the hard-core model in graph theory, especially via the local occupancy framework. Early on, in our discussion shortly following~\eqref{eq:occgreedy}, we mentioned that the condition we use cannot be improved in general. Put in another, perhaps blatantly obvious, way, the local occupancy condition is sharp for cliques in the most general (bounded-degree) circumstances, in which neighbourhoods are permitted maximum density.

We next discuss sharpness and barriers nearer the opposite end with respect to local density -- that is, triangle-free graphs. 
We first recap our analysis for triangle-free graphs. We established local $(\beta,\gamma)$-occupancy with specific functions $\beta$ and $\gamma$ of $\lambda$, whereby the optimization for a suitable choice of $\lambda$ tells us that $\beta+\gamma\Delta$ can be as small as $(1+o(1))\Delta/\log\Delta$ as $\Delta\to\infty$. As corollaries of Theorem~\ref{lem:occfrac}, Theorem~\ref{lem:fractional}, and Theorem~\ref{thm:colouring}, we deduced the following statements, for any triangle-free $G$ of maximum degree $\Delta$, in increasing order of strength, as $\Delta\to\infty$:
\begin{align*}
\alpha(G)\gtrsim \frac{n\log\Delta}{\Delta}, \qquad
\chi_f(G) \lesssim \frac{\Delta}{\log\Delta}, \qquad
\chi(G) \lesssim \frac{\Delta}{\log\Delta}.
\end{align*}

We next recall the off-diagonal Ramsey numbers and the asymptotic lower bound $R(3,k)\gtrsim k^2/(4\log k)$. We outlined how Shearer's bound on the independence number of triangle-free graphs of bounded maximum degree directly yields the best-known asymptotic upper bound, $R(3,k)\lesssim k^2/\log k$. As such, the local occupancy optimization for triangle-free graphs is sharp up to an asymptotic factor of {\em at most} $4$, and this is essentially due to the ultimate graph in the triangle-free process (see~\cite{Boh},~\cite{BoKe}, and~\cite{FGM}).
It is known that this (random) triangle-free graph has expected average degree $\Theta(\sqrt{n\log n})$.

The bounded-degree independence number bound, however, is more general than the Ramsey number bound, and it is matched more closely by other constructions in a regime with lower edge-density. We take $p = p(n)$ so that $d=np$ satisfies $d= \omega(\log n)$ and $d= o(n^{1/3})$, and consider the binomial random graph ${\mathbb G}_{n,p}$. Classic first moment estimates (see, for example,~\cite{AlSpbook}) show that $\Delta({\mathbb G}_{n,p})\lesssim d$, that $\alpha({\mathbb G}_{n,p})\lesssim (2n\log d)/d$, and that ${\mathbb G}_{n,p}$ has $o(n)$ triangles, with probability tending to $1$ as $n\to\infty$. Thus, there is some triangle-free induced subgraph $G=(V,E)$ of ${\mathbb G}_{n,p}$ of maximum degree $\Delta\sim d$ with $|V|\sim n$, for which $\alpha(G) \lesssim (2|V|\log\Delta)/\Delta$ with probability tending to $1$ as $n\to\infty$. This triangle-free construction shows how, over this range of densities, the triangle-free local occupancy optimization is sharp up to an asymptotic factor of at most $2$.

As an aside, this specific factor $2$ is perhaps not unrelated to an entrenched asymptotic factor $2$ in a longstanding open problem in random graph theory. For any fixed $\eps>0$, there is a  simple polynomial-time algorithm to find, with probability tending to $1$ as $n\to\infty$, an independent set of size $(1-\eps)(n\log d)/d$ in ${\mathbb G}_{n,p}$ as above. But in 1976, Karp~\cite{Kar76} essentially asked what happens if we replace the factor $1-\eps$ by $1+\eps$. (In fact, Karp asked this in the dense case $p=\frac12$, but the problem is well founded and is just as difficult in sparser regimes for $p=p(n)$.) The `polynomial-time' aspect of the problem is crucial, as there are non-constructive methods that succeed even if we replace the factor $1-\eps$ by $2-\eps$ (see~\cite{KaMc15} and~\cite{CoEf15} for further background). It is relevant to note that the local occupancy framework is `algorithmic' subject to efficient local sampling from the hard-core model~\cite{behemothalgorithmic}.
As such, the framework must contend with the `algorithmic barrier' for independent sets and colourings in random graphs.

In fact, it is nearly the most `tree-like' regime for triangle-free graphs in which we see the barriers in our framework even more clearly. Recall that the optimization for triangle-free graphs and Theorem~\ref{lem:occfrac} show that with $\lambda=\lambda(\Delta)=1/\log\Delta$ (say) the hard-core model at fugacity $\lambda$ on a triangle-free graph $G$ of maximum degree $\Delta$ has occupancy fraction satisfying
\begin{align*}
\alpha_G(\lambda)
\gtrsim \frac{\log\Delta}{\Delta}
\end{align*}
as $\Delta\to\infty$,
where $X$ is the random independent set.
As the occupancy fraction is monotone increasing in terms of $\lambda$, we maintain the same occupancy fraction guarantee when $X$ is a {\em uniformly} random independent set:
\begin{align*}
\alpha_G(1)
\ge (1+o_{\Delta}(1)) \frac{\log\Delta}{\Delta}.
\end{align*}
Now, for fixed $\Delta\ge 3$, consider the hard-core model at fugacity $\lambda$ on ${\mathbb G}_{n,\Delta}$, the random $\Delta$-regular graph.
By the work of Bhatnagar, Sly and Tetali~\cite{BST16}, its occupancy fraction is related to the unique translation-invariant hard-core measure on the infinite $\Delta$-regular tree, for all sufficiently large $\Delta$. From this, it follows that 
\begin{align*}
\alpha_{{\mathbb G}_{n,\Delta}}(1)
= (1+o_{\Delta}(1)) \frac{\log\Delta}{\Delta}
\end{align*}
 with probability tending to $1$ as $n\to\infty$.
 Moreover, it is a well-known fact that ${\mathbb G}_{n,\Delta}$ is triangle-free (and even of large fixed girth) with positive probability as $n\to\infty$.
Thus, the random $\Delta$-regular graph {\em conditioned on being triangle-free} shows how, in the regime of constant edge-density and at $\lambda=1$, the triangle-free local occupancy optimization is asymptotically sharp as $\Delta\to\infty$ (see~\cite{DJPR17independent} for more details on this discussion).

\medskip

At the risk of over-generalization, one might discern a pattern through the applications described in this chapter, and perhaps even in the discussion of sharpness constructions in this section. While the local occupancy method gives interesting results for independent sets and colourings over a wide array of settings, it has its limitations and gives better quality results when the problem under consideration is better characterized by local structure. 

To illustrate this, here is an important related problem where we would be delighted to make significant improvements through such methods, but have been unable to do so. 
Let $G$ be a {\em bipartite} graph of maximum degree $\Delta$. How large, as a function of $\Delta$, can the list chromatic number $\chi_\ell(G)$ of $G$ be?
The rub here is that, while $G$ is necessarily triangle-free and thus has excellent local structure, the methods that we have used in this chapter do not `see' (and so cannot take advantage of) the assumed global structure, the division of the vertex-set into two independent sets. Alon and Krivelevich~\cite{AlKr98} have boldly conjectured that the answer to the above question is at most $O(\log \Delta)$ -- a quantity that is drastically smaller than $\Theta(\Delta/\log\Delta)$ which we have in the triangle-free case. Unfortunately, `non-local' methods have not made much headway either (see~\cite{ACK21} and~\cite{BraPhDthesis} for recent progress).

More subtly, let us recall two of the applications described in Section~\ref{sec:lococcanalysis}.
The corollaries following Theorems~\ref{lem:densitylocallysparse} and~\ref{lem:densityCkfree}
each smoothly generalize the triangle-free case, and thus skirt close to the barriers, but these bounds deteriorate as they stray from the triangle-free case. More precisely, they become no better than the trivial bounds if, respectively, we forbid having $t$ edges induced in any neighbourhood for some $t=\Omega(\Delta^2)$, or we forbid as a subgraph a cycle $C_k$ of length $k=\Omega(\Delta)$. Thus, the local occupancy method seems to break down with non-vanishing local edge-density. On the other hand, we still know decent bounds in this regime -- that is, there is a proper colouring using at most $(1-\eps)\Delta$ colours for some positive fixed $\eps$, where $\eps$ depends on the local edge-density bound. It might come as a surprise that a seemingly blunter method succeeds in this case: the so-called `na\"ive colouring' method, whereby each vertex is assigned a uniformly random colour from some ground set of colours and then conflicts (that is, monochromatic edges) are resolved by a simple follow-up recolouring procedure. See~\cite{molloy_graph_2002} for a broader treatment of this method, and~\cite{HDK22} for the most recent developments along such lines.

\section{Open problems}

We conclude by highlighting some open problems near the limits of local occupancy.

The first problem concerns occupancy fraction.
Given $\lambda>0$, let $f_\lambda(d)=\lambda/(1+(d+1)\lambda)$. We showed in Section~\ref{sec:intro} that, given a graph $G$ of maximum degree $\Delta$, the occupancy fraction of the hard-core model at fugacity $\lambda$ satisfies 
$\alpha_G(\lambda) \ge f_\lambda(\Delta)$,
and that this is sharp when $G$ is a disjoint union of $K_{\Delta+1}$s. We propose the following conjecture.
\begin{conjecture}\label{conj:hcmcarowei}
With $f_\lambda$ as above, let $X$ be a random independent set in a graph $G=(V,E)$, chosen according to the hard-core model at fugacity $\lambda$. Then
\begin{align*}
\EE |X| \ge \sum_{v\in V} f_\lambda(\deg(v)).
\end{align*}
\end{conjecture}
\noindent
As $f_\lambda$ is convex, this summation is bounded from below by $f_\lambda(\overline{\deg}(G))|V|$, where $\overline{\deg}(G)$ is the average degree of $G$. Thus, Conjecture~\ref{conj:hcmcarowei} is a stronger form of the first occupancy fraction bound we gave, while remaining sharp when $G$ is a disjoint union of complete graphs. Its proof would moreover simultaneously refine and strengthen the bound in~\eqref{eq:EXlbavgdeg}, the Caro--Wei theorem (see~\cite{Car79} and~\cite{Wei81}), and an analogous lower bound on $\log Z_G(\lambda)$ given by Sah, Sawhney, Stoner and Zhao~\cite[Theorem~1.7]{sah_number_2019}. Conjecture~\ref{conj:hcmcarowei} implies the latter result by integration (see~\eqref{eq:Z'}).

The second problem concerns the relationship between occupancy fraction and the independence number. We have already observed how, in the complete graph $K_{\Delta+1}$, the average size of an independent set is $(\Delta+1)/(\Delta+2)$, while the independence number is $1$, and so their ratio approaches $1$ as $\Delta\to\infty$.
Because of the barriers discussed above, it is impossible to improve upon Shearer's bound $R(3,k)\lesssim k^2/\log k$ by a na\"ive application of Theorem~\ref{thm:tflococc}, which requires small $\lambda$ for the best quantitative bounds. 
But there remains a route to move further beyond local occupancy 
-- that is, if, for all triangle-free graphs $G$, the ratio between $\alpha(G)$ and the average size of an independent set in $G$
is bounded above $1$, then one obtains an improvement to the constant factor in Shearer's bound. 
The following conjecture from~\cite{DJPR17average} asserts that, in triangle-free graphs, the ratio approaches $2$ as the minimum degree grows. 
By a simple argument that removes small-degree vertices, the following conjecture if true would imply that $R(3,k)\lesssim k^2/(2\log k)$.

\begin{conjecture}
\label{conj:tfavgmax}
If $G$ is a triangle-free graph of minimum degree $d$, then  as $d\to\infty$
\begin{align*}
\frac{\alpha(G)}{\alpha_G(1)|V(G)|} 
\ge 2-o_d(1).
\end{align*}
\end{conjecture}

Last, let $G$ be a $K_4$-free graph of maximum degree $\Delta$, so that every neighbourhood of $G$ induces a triangle-free graph. Under this constraint, there may still be non-trivial edge-density in neighbourhoods, and this falls into more difficult territory; however, in any case a suitable local occupancy optimization (see~\cite{davies_graph_2020}) yields an independent set of size $(1+o_\Delta(1))(n\log \Delta)/(\Delta\log\log\Delta)$, which is currently the best bound. There is some distance to cover, as the longstanding conjecture of Ajtai, Erd\H{o}s, Koml\'os and Szemer\'edi~\cite{AEKS81} posits much more.

\begin{conjecture}
\label{conj:AEKS81}
Given a $K_4$-free graph $G$ of maximum degree $\Delta$ on $n$ vertices, there is an independent set of size at least $(c n\log \Delta)/\Delta$ for some absolute constant $c>0$.
\end{conjecture}

\paragraph{Acknowledgements} Ewan Davies was supported in part by NSF grant CCF-2309707. Ross Kang was supported by Vidi (639.032.614) and Open Competition (OCENW.M20.009) grants of the Dutch Research Council (NWO) as well as  the Gravitation programme NETWORKS (024.002.003) of the Dutch Ministry of Education, Culture and Science (OCW). 

\paragraph{Notes added}
We posted a version of this book chapter in a public preprint repository in January of 2025. Not long afterwards came a number of interesting related developments.

In March 2025 Buys, van den Heuvel and Kang~\cite{BHK25} generalized Shearer's proof~\cite{She83} that any $n$-vertex triangle-free graph of {\em average} degree $d$ contains an independent set of size $(1-o(1))n\log(d)/d$ to the partition function of the hard-core model. 
This yields the lower bound that the number of independent sets (of all sizes) in such a graph must be at least 
\[ \exp\left((1-o(1))\frac{(\log d)^2}{2d}n\right),\]
which is best possible (keeping in mind~\eqref{eq:Z'}).

In April 2025 Klartag~\cite{Kla25} defined a randomized construction of a large origin-symmetric ellipsoid in $\mathbb{R}^d$ containing no points of $\mathbb{Z}^d$ other than the origin.
His analysis of the volume of this ellipsoid gives the existence of lattice-based sphere packings that certify the new density lower bound $\theta(d)\ge\Omega(d^2\cdot 2^{-d})$, improving on the result in~\cite{campos_new_2023}.

While Shearer's remains the state of the art for  upper bounds on $R(3,k)$, in May 2025 Campos, Jenssen, Michelen and Sahasrabudhe~\cite{CJMS25} improved the lower bound on $R(3,k)$ to $(1/3-o(1))k^2/\log k$ as $k\to\infty$ by analyzing a variant of the triangle-free process that starts from a somewhat structured graph. 
It seems that starting the process carefully and steering it along a simpler trajectory is easier to analyze and yields a better bound. 

Also in May 2025, Davies, Sandhu and Tan~\cite{DST25} showed that bounds of the form~\eqref{eq:lococc} can be used to make progress on Conjecture A. They established the conjecture in graphs of maximum degree $\Delta$ for values $\lambda\le 3/(\Delta+1)^2$ (alternatively, for each $n\ge 1$ the conjecture holds in all $n$-vertex graphs for $\lambda\le 3/n^2$).


\footnotesize
\bibliographystyle{robin}
\bibliography{references}

\end{document}